\documentclass[authoryear]{elsarticle}
\usepackage[utf8]{inputenc}

\usepackage{subfig}
\usepackage{amsmath}
\usepackage{amsfonts}
\usepackage{multirow}
\usepackage{tabularx}
\usepackage{booktabs}
\usepackage{xltabular}
\usepackage{appendix,setspace}

\usepackage{lscape}

\usepackage{listings}
\lstdefinestyle{ExStyle}{
moredelim=**[is][\color{red}]{~}{~},
}

\usepackage{tikz}
\usetikzlibrary{calc, shapes, snakes, arrows, arrows.meta, positioning}
\usepackage{pgfplots}
\pgfplotsset{compat=1.13}

\newcommand*\galini{\mbox{GALINI}}

\newcommand{\xil}{y_{il}}
\newcommand{\ylj}{y_{lj}}
\newcommand{\yij}{y_{ij}}
\newcommand{\xij}{y_{ij}}

\newcommand{\plk}{p_{lk}}
\newcommand{\qil}{q_{il}}
\newcommand{\ci}{c_i}
\newcommand{\rj}{d_j}
\newcommand{\Ail}{A_i^L}
\newcommand{\Aiu}{A_i^U}
\newcommand{\Sl}{S_l}
\newcommand{\Djl}{D_j^L}
\newcommand{\Dju}{D_j^U}
\newcommand{\Cik}{C_{ik}}

\newcommand{\Pjkl}{P_{jk}^L}
\newcommand{\Pjku}{P_{jk}^U}
\newcommand{\clj}{c_{lj}}

\newcommand{\zlj}{z_{lj}}
\newcommand{\slj}{s_{lj}}
\newcommand{\uljk}{u_{ljk}}
\newcommand{\tljk}{t_{ljk}}
\newcommand{\pljk}{p_{ljk}}
\newcommand{\rljk}{r_{ljk}}
\newcommand{\betalowerljk}{\underline{\beta}_{ljk}}
\newcommand{\betaupperljk}{\overline{\beta}_{ljk}}
\newcommand{\etalowerljk}{\underline{\eta}_{ljk}}
\newcommand{\etaupperljk}{\overline{\eta}_{ljk}}

\newcommand{\wiltj}{w_{il_tj}}
\newcommand{\zltj}{\zeta_{l_tj}}
\newcommand{\vilj}{v_{ilj}}
\newcommand{\gammalt}{\gamma_{lt}}

\begin{document}

\onehalfspacing

\begin{frontmatter}

\title{Solving the pooling problem at scale with extensible solver \galini}

\author[icl]{F. Ceccon}
\author[icl]{R. Misener\corref{cor1}}
\address[icl]{Department of Computing, Imperial College London, 180 Queens Gate, SW7 2AZ, UK}

\cortext[cor1]{\texttt{r.misener@imperial.ac.uk}, +44 (0) 20759 48315}

\begin{abstract}
  This paper presents a Python library to model pooling problems, a class of network flow problems with many engineering applications.
  The library automatically generates a mixed-integer quadratically-constrained quadratic optimization problem from a given network structure.
  The library additionally uses the network structure to build 1) a convex linear relaxation of the non-convex quadratic program and 2) a mixed-integer linear restriction of the problem.
  We integrate the pooling network library with \galini{}, an open-source extensible global solver for quadratic optimization. We demonstrate \galini's extensible characteristics by using the pooling library to develop two \galini{} plug-ins: 1) a cut generator plug-in that adds valid inequalities in the \galini{} cut loop and 2) a primal heuristic plug-in that uses the mixed-integer linear restriction.
  We test \galini{} on large scale pooling problems and show that, thanks to the good upper bound provided by the mixed-integer linear restriction and the good lower bounds provided by the convex relaxation, we obtain optimality gaps that are competitive with Gurobi 9.1 on the largest problem instances.
\end{abstract}


\end{frontmatter}

\section{Introduction}

This manuscript presents a library to formulate pooling problems using Pyomo~\citep{Hart2017,Watson2011}, a modeling library for Python, without losing network structure information.
Traditionally, mixed-integer non-linear optimization (MINLP) solvers see ``flat'' optimization problems as a series of objectives, constraints, and variables, with little to no information about what the optimization problem represents.
We propose retaining the network structure: we use Python to model the network structure with nodes representing network components, and edges representing possible network flow.

We focus on the pooling problem because of its many industrial applications, including \citep{misener2009advances}: crude-oil scheduling \citep{lee1996mixed,li2007improving,li2012continuous,li2012scheduling}, water networks \citep{galan1998optimal,castro2007efficient}, natural gas production \citep{selot2008short,li2011stochastic}, fixed-charge transportation with product blending \citep{papageorgiou2012fixed}, hybrid energy systems \citep{baliban2012global}, multi-period blend scheduling \citep{kolodziej2013discretization}, and mining \citep{boland2015special}. 
Solving the pooling problem is NP-hard \citep{alfaki2013strong,Baltean-Lugojan2018,letsios2020approximation}, so deterministic global optimization solvers use algorithms such as branch \& bound to solve the problem. Our goal with explicitly using special structure information is practically solving larger problem sizes.
But embedding structural information within \galini{} could be similarly implemented for other optimization problems over networks, e.g.\ optimal power flow.

Our pooling network library automatically generates a Pyomo block with the problem formulated as a mixed-integer quadratically constrained program (MIQCQP). The library also provides functions to compute a feasible (but not necessarily optimal) problem solution \citep{Dey2015}, and functions to generate valid linear cuts \citep{luedtke2020strong}.
We integrate the pooling network library with \galini{} \citep{ceccon2020galini, ceccon2021galini}, an extensible MIQCQP solver. The pooling network library provides an initial feasible solution search strategy and a cut generator: users enable these components by changing their \galini{} configuration file.
Both these \galini{} extensions take advantage of the pooling problem embedded inside the Pyomo model. \galini{} itself does not require any change, so everything presented in this paper is an external extension to \galini{}. The advantage of this approach is that, although \galini{}  itself is unaware of the pooling problem, \galini{} can, through this library, behave as if it has perfect knowledge of its structure.

State-of-the-art MIP solvers detect constraints that commonly appear in network problems and deduce the network structure \citep{brown1984automatic,bixby1988finding,gulpinar2004extracting,achterberg2010mcf,salvagnin2016detecting}.
Out prior work \citep{Ceccon2016} deduced MINLP network structure: we detected pooling problems in 6\% of MINLPLib. But detecting special structure in a ``flat'' optimization problem is a heuristic approach that may miss special structure developed by the user.
Our pooling network library is different because the library keeps the original pooling problem structure together with the optimization problem. 
This paper is published alongside the source code \citep{ceccon2021pooling}, licensed under the Apache License, version 2.

\section{Background}

\subsection{The pooling problem}

The pooling problem has a feed-forward network $T = (N, A)$, where $N$ and $A$ are the sets of nodes and edges, respectively. Nodes $N$ are partitioned in three layers: 1) inputs $I$ representing material availability, 2) pools $L$ representing material mixing, and 3) outputs $J$ representing product demand. By definition, $N = I \cup L \cup J$. The network allows connections from inputs to pools, inputs to outputs, and pools to outputs, i.e.\ $A = (I \times L) \cup (I \times J) \cup (L \times J)$. The notation $I_l$ represents the input nodes connected to pool $l$, i.e.\ $I_l = \{i : i \in I, (i, l) \in A\}$. Similarly, $I_j$ represents the set of inputs connected to output $j$.
The optimization objective is to minimize cost (maximize profit) in the network, while tracking material qualities $K$. Tracking these material qualities K across the network introduces bilinear terms to the problem.

The first optimization formulation to solve the pooling problem was the P-formulation \citep{haverly1978studies}.
The mathematically-equivalent Q-formulation \citep{ben1994global} introduces fractional flow variables $\xil = \qil \sum_{j \in J} \ylj$.
The Equation~\eqref{eq:pooling.pq.formulation} PQ-formulation \citep{quesada1995global,tawarmalani2013convexification} is equivalent to the Q-formulation with the addition of Equation~\eqref{eq:pq.cut}. Tables~\ref{table:pooling.problem.notation.sets} and \ref{table:pooling.problem.notation.variables} in \ref{app:notation} summarize the notation.
\allowdisplaybreaks{
\begin{subequations}
  \label{eq:pooling.pq.formulation}
  \begin{align}
    & & & \underset{y,v,q}{\text{min}} \sum_{i \in I, l \in L} c_i \vilj - \sum_{l \in L, j \in J} d_j \ylj - \sum_{i \in I, j \in J} (d_j - c_i) \yij \\
    & &\begin{array}{r} \text{Path} \\ \text{Definition} \end{array} &\left[ v_{ilj} = q_{il} y_{lj} \quad \forall i \in I, l \in L, j \in J
                                                                       \vphantom{\sum_{i \in I}} \right. \label{eq:pooling.pq.pathdef}\\
    & & \begin{array}{r} \text{Simplex} \end{array} & \left[ \sum_{i \in I} q_{il} = 1 \quad \forall l \in L \right. \label{eq:pooling.pq.simplex}\\
    & & \begin{array}{r} \text{Product} \\ \text{Quality} \end{array}&
                                                                       \left[
                                                                       \begin{array}{lr}
                                                                       \displaystyle \sum_{i \in I, l \in L} \Cik \vilj + \displaystyle \sum_{i \in I} \Cik \yij  & \left\{
                                                                       \begin{array}{l}
                                                                         \leq \Pjku \left( \displaystyle \sum_{l \in L} \ylj + \sum_{i \in I} \yij \right)\\
                                                                         \geq \Pjkl \left( \displaystyle  \sum_{l \in L} \ylj + \sum_{i \in I} \yij \right)
                                                                       \end{array}
    \right.\\
    & \quad \forall k \in K, j \in J
    \end{array}
     \right. \label{eq:pooling.pq.product.quality}\\
    & & \begin{array}{r} \text{Input} \\ \text{Capacity} \end{array}&\left[ A_i^L \leq \sum_{l \in L, j \in J} \vilj + \sum_{j \in J} \yij \leq A_i^U \quad \forall i \in I \right.\label{eq:pooling.pq.input.capacity}\\
    & & \begin{array}{r} \text{Pool} \\ \text{Capacity} \end{array}&\left[ \sum_{l \in L} \ylj \leq S_l \quad \forall l \in L\right. \label{eq:pooling.pq.pool.capacity}\\
    & & \begin{array}{r} \text{Output} \\ \text{Capacity} \end{array}&\left[ D_j^L \leq \sum_{l \in L} \ylj + \sum_{i \in I} \yij \leq D_j^U \quad \forall j \in J \right. \label{eq:pooling.pq.output.capacity}\\
    & & \begin{array}{r} \text{Reduction} \end{array}&\left[\sum_{i \in I} v_{ilj} = y_{lj} \quad \forall l \in L, j \in J \right. \label{eq:pooling.pq.reduction.1}\\
    & & \begin{array}{r} \text{Reduction} \end{array}&\left[\sum_{i \in I, l \in L} \vilj \leq  c_l \qil \quad \forall i \in I, l \in L \right. \label{eq:pooling.pq.reduction.2}\\
    & & \begin{array}{r} \text{Other} \\ \text{Bounds} \end{array}&\left[
                                                                    \begin{array}{l}
                                                                      \sum_{j \in J} v_{ilj} \leq c_{il} \quad \forall i \in I, l \in L\\
                                                                      y_{lj} \leq c_{lj} \quad \forall l \in L, j \in J\\
                                                                      y_{ij} \leq c_{ij} \quad \forall i \in I, j \in J\\
                                                                      y_{ij}, y_{lj}, v_{ilj}, q_{il} \geq 0 \quad \forall i \in I, l \in L, j \in J\\
                                                                      \sum_{j \in J} v_{ilj} \leq c_l q_{il} \quad \forall i \in I, l \in L
                                                                    \end{array}
    \right.
    \label{eq:pooling.pq.bounds}
  \end{align}
\end{subequations}
}

\begin{equation}
  \label{eq:pq.cut}
  \sum_{i \in I} \qil \ylj = \ylj \quad \forall l \in L, j \in J
\end{equation}

Other pooling formulations incorporate other algorithmic trade-offs \citep{audet2004pooling,alfaki2013multi,boland2016new}.

\subsection{Pyomo}

Pyomo \citep{Hart2017, Friedman2013} is a Python-based algebraic modeling language to formulate MINLP. Other algebraic modeling languages include GAMS \citep{brooke1996gams}, AMPL \citep{fourer2003ampl}, AIMMS, YALMIP \citep{lofberg2004yalmip}, and JuMP \citep{dunning2017jump}.

Large scale optimization models describe complex physical systems comprised of many sub-systems.
Pyomo provides a \texttt{Block} component to structure models. Blocks encapsulate individual model components in reusable pieces that are combined together to form the final model. For example, users can develop blocks representing their system at each time period and combine multiple instances of the block to produce a scheduling problem.

\begin{figure}
  \centering
  \includegraphics[width=1.0\textwidth]{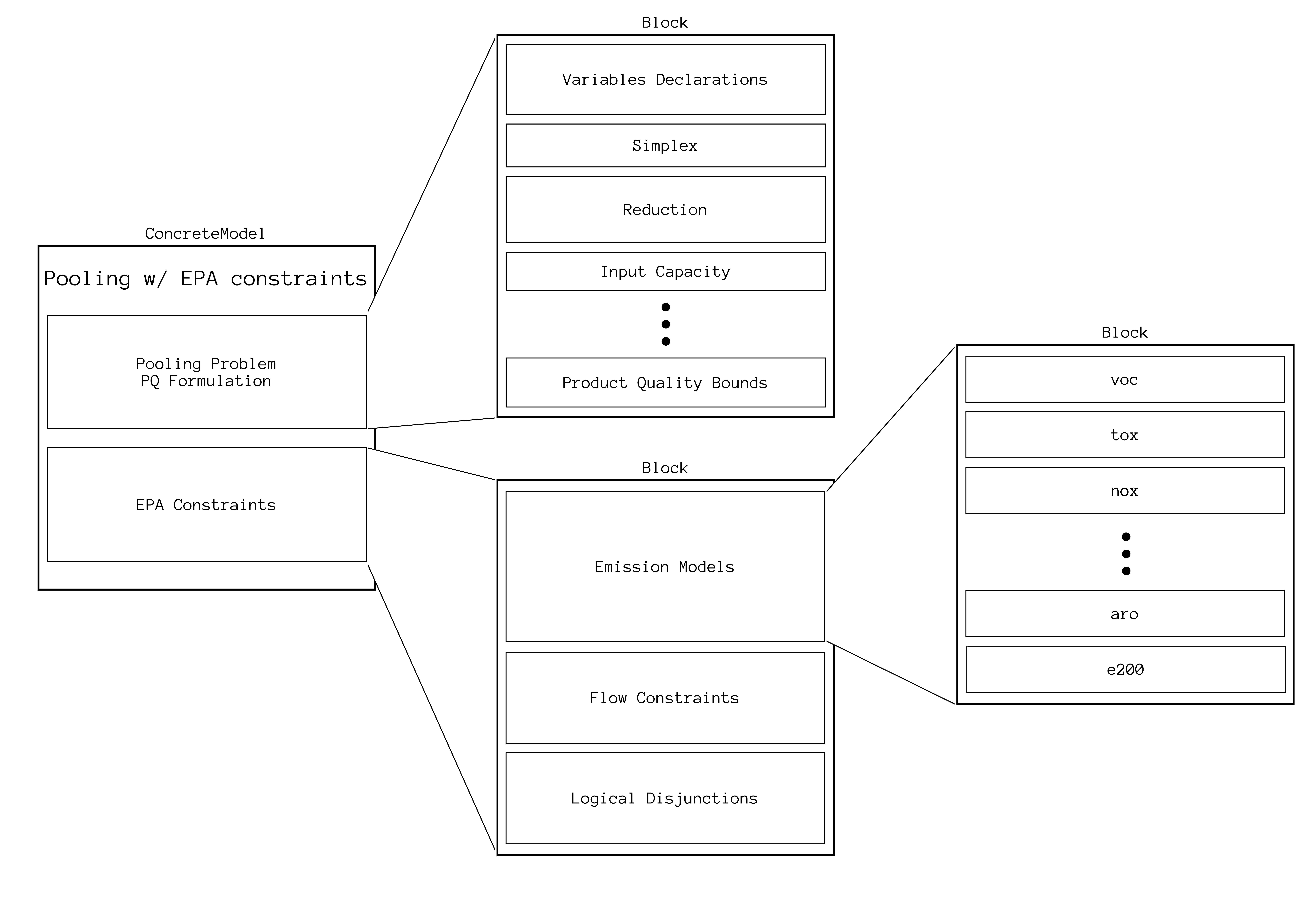}
  \caption{Pyomo model for the pooling problem with EPA constraints \citep{Misener2010d}. The model contains one reusable model with the pooling problem PQ-formulation and one block with the additional EPA constraints. Each block contains definitions for its variables and constraints, the EPA block uses variables from the PQ-formulation block.}
  \label{fig:pyomo.block.pooling}
\end{figure}

Pyomo blocks can be nested into other blocks, creating a tree-like structure of blocks. Figure~\ref{fig:pyomo.block.pooling} shows how to organize a large model to optimize refinery production subject to US Environmental Protection Agency (EPA) regulation \citep{Misener2010d}. The root \texttt{ConcreteModel} has two blocks, one containing the refinery production model, the other containing EPA constraints. The EPA constraints contain complex non-linear expressions and logical disjunctions, so we organize them further in separate blocks. We add one block for each emission model to make the model easier to update and debug. Blocks may use variables defined in other blocks, both from their children blocks and their siblings.

Developers can also create \textit{smart} blocks, i.e.\ new types of Pyomo blocks, which can create and update their internal constraints. These custom blocks can automatically generate relaxations of a particular expression or build an optimization model from an higher-level definition. \ref{app:pyomo_smart_block} shows how to define a new Pyomo block that relaxes a nonconvex bilinear term.

\subsection{\galini{}}

\galini{} is an open-source MIQCQP solver written in Python that can be extended at runtime.
\galini{} implements a generic branch \& cut algorithm: at each node of the branch \& bound algorithm \galini{} iteratively adds cuts and solves the linear relaxation of the original optimization problem.
Algorithm developers extend \galini{} by registering their classes using Python \textit{entry points}, a system-wide registry of Python classes.
This manuscript extends \galini{} with a new cut generator and an initial primal search heuristic. \galini{} uses cut generators inside the cut loop to augment the linear relaxation of the optimization problem. \galini{} calls the initial primal search heuristic before visiting the branch \& bound root node to find an initial feasible solution.

The cut generator interface requires to implement six methods:

\begin{itemize}
\item \texttt{before\_start\_at\_root(problem, relaxed\_problem)}: called before entering the cut loop at the root node,
\item \texttt{after\_end\_at\_root(problem, relaxed\_problem, solution)}: called after the cut loop at the root node,
\item \texttt{before\_start\_at\_node(problem, relaxed\_problem)}: called before entering the cut loop  at non-root nodes,
\item \texttt{after\_end\_at\_node(problem, relaxed\_problem, solution)}: called after the cut loop at non-root nodes,
\item \texttt{has\_converged(state)}: returns a true value if the cut generator won't generate any more cuts at this node,
\item \texttt{generate(problem, relaxed\_problem, solution, tree, node)}:
  return a list of cuts to be added to \texttt{relaxed\_problem}, can use information from the previous iteration \texttt{solution}.
\end{itemize}

The \texttt{InitialPrimalSearchStrategy} requires algorithm developers to implement the following method:

\begin{enumerate}
\item \texttt{\textbf{solve(model, tree, node)}}: find a feasible solution for the model if any, or \texttt{None} if not possible.
\end{enumerate}

\section{Pooling network library}
\label{chapt:network.structure}

\subsection{Network structure}

\begin{figure}%
  \centering
  \includegraphics[width=0.90\textwidth]{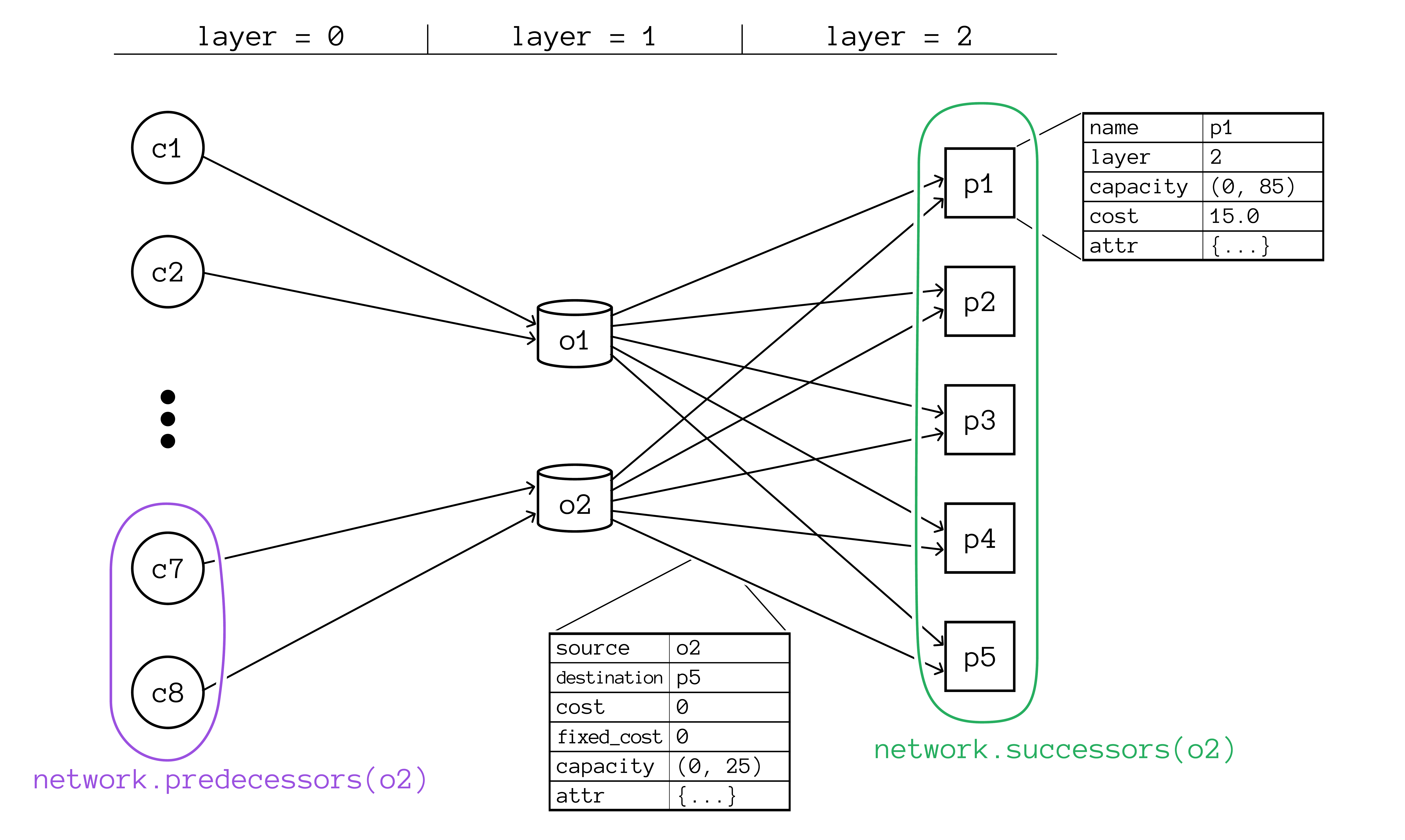}
  \caption{Graphical representation of the pooling problem \texttt{adhya4} as implemented with the pooling network library. Nodes are grouped by layer, in this example there are 3 layers (one for inputs, one for pools, and one for outputs). Nodes have the attributes: \texttt{name}, \texttt{layer}, \texttt{capacity}, and \texttt{cost}, and the (optional) \texttt{attr} dictionary to add application-specific attributes. Edges have the attributes: \texttt{source}, \texttt{destination}, \texttt{cost} (per unit flow cost), \texttt{fixed\_cost}, \texttt{capacity}, and the (optional) \texttt{attr} dictionary to add application-specific attributes The figure shows the predecessors and successors of node \texttt{o2} highlighted in purple and green, respectively.}%
  \label{fig:pooling.problem.network.library}%
\end{figure}

The library's \texttt{Network} class describes the network structure of the problem. This manuscript assumes that a user explicitly stores the network structure. But we could alternatively detect pooling network structure \citep{Ceccon2016} and then use pooling network library to develop algorithmic approaches, e.g. cutting planes and primal heuristic.
\ref{app:network_structure} shows how the pooling library builds and stores a problem. As an example, Figure~\ref{fig:pooling.problem.network.library} shows how the library represents the \texttt{adhya4} test problem \citep{adhya1999lagrangian}. 


\subsection{Automatic PQ-formulation generation}

\begin{figure}%
  \centering
  \includegraphics[width=0.40\textwidth]{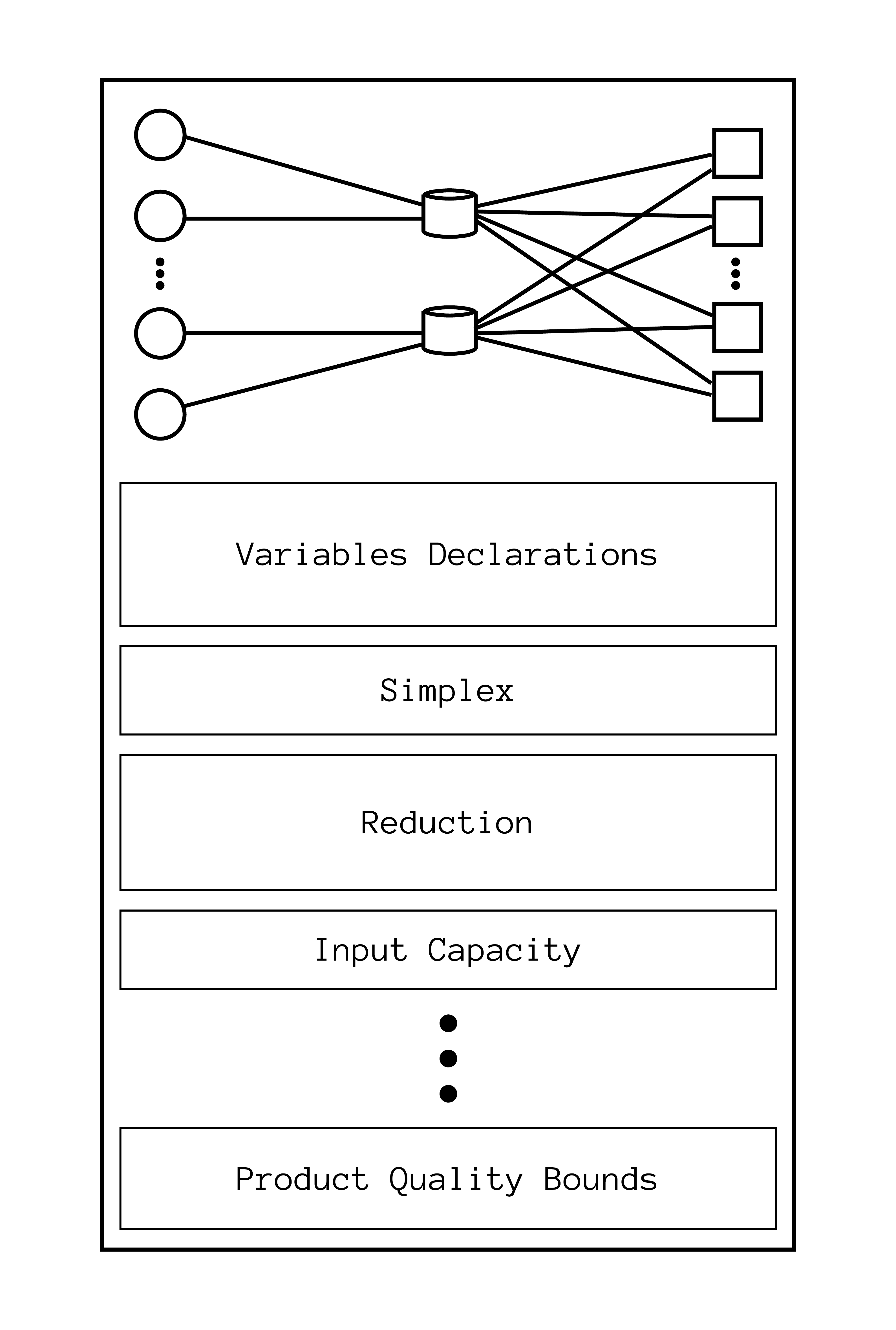}
  \caption{\texttt{PoolingPQFormulation} Pyomo block. The block automatically generates the variables and constraints for the pooling PQ-formulation. The block keeps a reference to the pooling problem network structure. Software that is aware of this structure can use it to have a precise knowledge of the pooling problem.}%
  \label{fig:pooling.problem.network.block}%
\end{figure}

The library provides a custom \texttt{PoolingPQFormulationBlock} Pyomo block to generate the pooling problem PQ-formulation in Optimization Problem~\eqref{eq:pooling.pq.formulation} from the network structure. The library does not check that the network passed to it is a valid pooling problem, but it expects a network with a pooling problem-like structure. If the network contains edges that are not part of a standard pooling problem, the library ignores these edges. The network is expected to have a three layer structure. If the network has more than three layers of nodes, then all the layers after the third are ignored.
The pooling PQ-formulation block can be used in an indexed Pyomo block to have multiple pooling problems in the same optimization problem. The block adds the variables and constraints to the block, so, after changing anything in the optimization problem, the PQ-formulation needs to be rebuilt by calling the \texttt{rebuild()} method on the block. Table~\ref{table:block.formulation} summarizes the block's variables and how they relate to the PQ-formulation. 
Figure~\ref{fig:pooling.problem.network.block} shows the structure of a \texttt{PoolingPQFormulationBlock} Pyomo block, with the pooling problem constraints stored together with the network structure.

\begin{table}
  \begin{center}
    \begin{tabularx}{\textwidth}{@{} l l X @{}}
      \toprule
      Type & Name & Description \\
      \midrule
      Variables & \texttt{q[i, l]} & Fractional flow: input $i$ to pool $l$, as fraction of total flow through $l$ \\
      & \texttt{v[i, l, j]} & Flow: input $i$ to output $j$ through pool $l$ \\
      & \texttt{y[l, j]} & Pool $l$ to output $j$ flow \\
      & \texttt{z[i, j]} & Input $i$ to output $j$ flow \\
      \midrule
      Constraints & \texttt{path\_definition[i, l, j]} & Equation~\eqref{eq:pooling.pq.pathdef}\\
      & \texttt{simplex[l]} & Equation~\eqref{eq:pooling.pq.simplex}\\
      & \texttt{product\_quality\_lower\_bound[j, k]} & Equation~\eqref{eq:pooling.pq.product.quality} \\
      & \texttt{product\_quality\_upper\_bound[j, k]} & Equation~\eqref{eq:pooling.pq.product.quality} \\
      & \texttt{input\_capacity[i]} & Equation~\eqref{eq:pooling.pq.input.capacity} \\
      & \texttt{pool\_capacity[l]} & Equation~\eqref{eq:pooling.pq.pool.capacity} \\
      & \texttt{output\_capacity[j]} & Equation~\eqref{eq:pooling.pq.output.capacity} \\
      & \texttt{reduction\_1[l, j]} & Equation~\eqref{eq:pooling.pq.reduction.1}\\
      & \texttt{reduction\_2[i, l]} & Equation~\eqref{eq:pooling.pq.reduction.2}\\
      \bottomrule
    \end{tabularx}
  \end{center}
  \caption{Top-level definitions of \texttt{PoolingPQFormulationBlock}. Users can use the block variables in constraints outside of the block. Users can also \texttt{deactivate()} the constraints.}
  \label{table:block.formulation}
\end{table}

The library also provides convenience functions to iterate over the indexes used in a pooling problem. For example, the function \texttt{index\_set\_ilj} will return an iterator over the indexes \texttt{(i, l, j)} of variables \texttt{v[i, l, j]}, equivalent to the variables $v_{ilj}$ in the pooling problem PQ-formulation.

\subsection{Convex relaxation and cut generator}

The pooling network library implements a convex relaxation of the pooling problem \citep{luedtke2020strong}.
The relaxation introduces a variable $\zlj$ that represents the total flow to output $j$ that does not pass through pool $l$:
\begin{equation*}
  \label{eq:pooling.cut.z.def}
  \zlj := \sum_{i \in I} \xij + \sum_{l' \in L, l' \neq l} y_{l'j} \quad \forall l \in L, j \in J
\end{equation*}
\noindent
The contribution of flow through pool $l$ to the quality $k$ at output $j$ is defined:
\begin{equation*}
  \label{eq:poolin.cut.u.def}
  \uljk := \sum_{i \in I} (\Pjku - \Cik)\vilj \quad \forall l \in L, j \in J, k \in K
\end{equation*}
\noindent
The total flow not going through pool $l$ to the quality $k$ at output $j$ is:

\begin{equation*}
  \label{eq:pooling.cut.t.def}
  \tljk := \sum_{i \in I} (\Pjku - \Cik) \xil + \sum_{i \in I} \sum_{l' \in L, l' \neq l} (\Pjku - \Cik) v_{il'j}
  \quad \forall l \in L, j \in J, k \in K
\end{equation*}
\noindent
The quality of attribute $k$ at pool $l$ is:
\begin{equation*}
  \label{eq:pooling.cut.p.def}
  \pljk := \sum_{i \in I} (\Pjku - \Cik) \qil \quad \forall l \in L, j \in J, k \in K
\end{equation*}
\noindent
Scaled flow variable $\slj$ scales the flow from pool $l$ to output $j$ to be in $[0, 1]$:
\begin{equation*}
  \label{eq:pooling.cut.s.def}
  \slj := \frac{1}{c_j}\sum_{i \in I} \vilj \quad \forall l \in L, j \in J
\end{equation*}
\noindent
We add the McCormick envelope of $\rljk = \slj \pljk$, where $\pljk^L$ and $\pljk^U$ are the bounds on $\pljk$:
\begin{align}
  \label{eq:pooling.cut.uxp}
  \rljk - \pljk^L \slj & \geq 0 \quad \forall l \in L, j \in J\\
  \pljk^U \slj - \rljk & \geq 0 \quad \forall l \in L, j \in J\\
  \rljk - \pljk^L x & \leq \pljk - \pljk^L \quad \forall l \in L, j \in J\\
  \pljk^U \slj  - \rljk & \leq \pljk^U - \pljk \quad \forall l \in L, j \in J
\end{align}
Parameters $\etalowerljk$ and $\etaupperljk$ track the lower and upper bounds on the excess of quality $k$ at output $j$ over the inputs $i$ that are connected to pool $l$.
\begin{align}
  \etalowerljk & := \min\{(\Pjku - \Cik) : i \in I_l\}  \quad \forall l \in L, j \in J\\
  \etaupperljk & := \max\{(\Pjku - \Cik) : i \in I_l\} \quad \forall l \in L, j \in J
\end{align}
\noindent
Similarly, the parameters $\betalowerljk$ and $\betaupperljk$ track the bounds on the excess of quality $k$ at output $j$ over the inputs $i$ that are \textit{not} connected to pool $l$.
\begin{align}
  \betalowerljk & := \min\{(\Pjku - \Cik) : i \in I_j \cup \bigcup_{l' \in L \setminus \{l\}} I_{l'}\} \quad \forall l \in L, j \in J\\
  \betaupperljk & := \max\{(\Pjku - \Cik) : i \in I_j \cup \bigcup_{l' \in L \setminus \{l\}} I_{l'}\} \quad \forall l \in L, j \in J
\end{align}

\begin{figure}%
  \centering
  \includegraphics[width=0.6\textwidth]{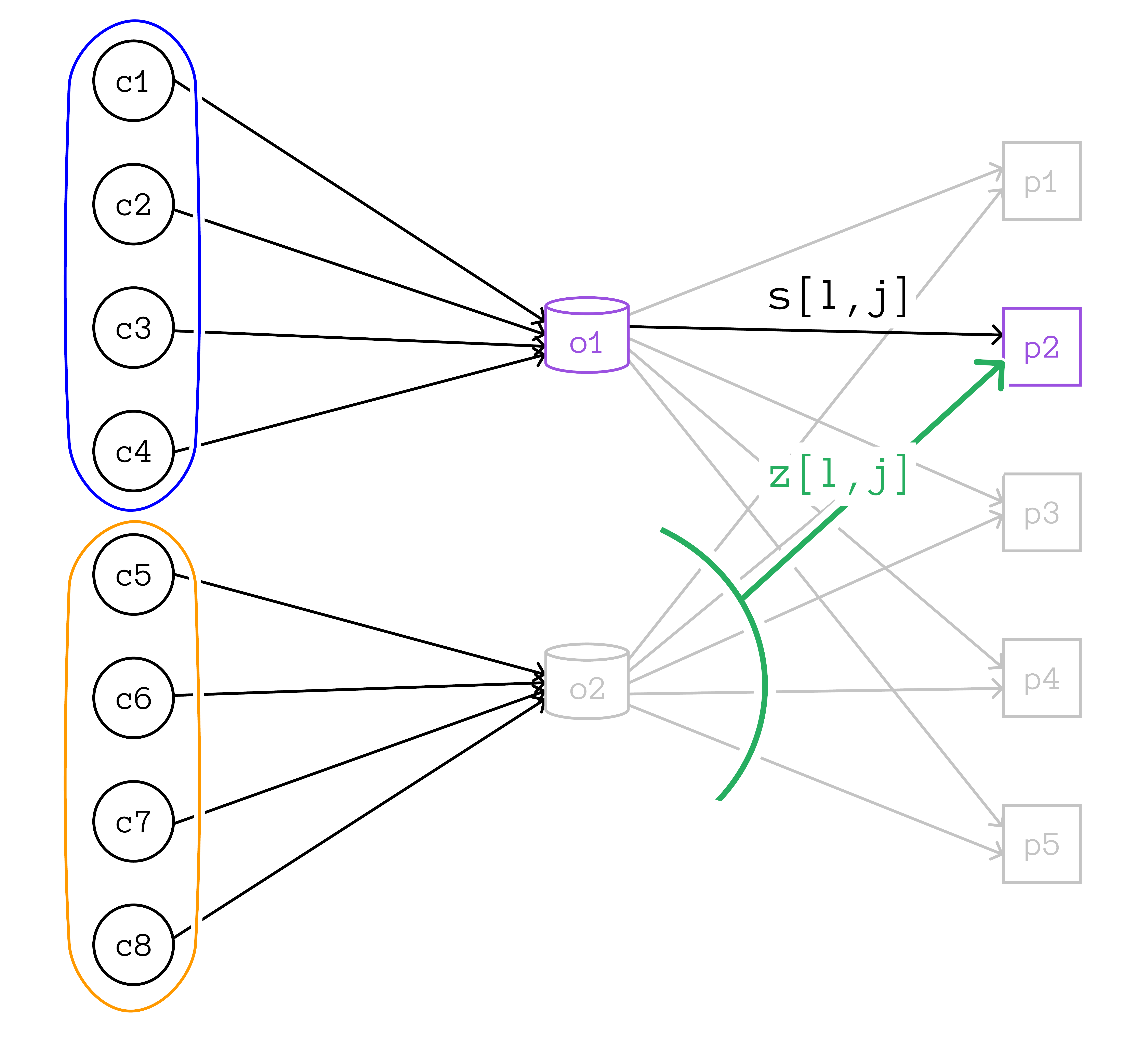}
  \caption{Variables and sets generating cuts for each pool, output, quality triplet. This example considers pool \texttt{o1} and output \texttt{p2}. Variable \texttt{s[l, j]} is the scaled flow from pool $l$ to output $j$. Variable \texttt{z[l, j]} is the flow to output $j$ not passing through pool $l$, i.e.\ the flow to output \texttt{p2} passing through pool \texttt{o2}. The inputs circled in blue are the inputs connected to the output $j$ through pool $l$ that compute $\etalowerljk$ and $\etaupperljk$. The inputs circled in orange are the inputs connected to output $j$ through the by-pass flow, and compute $\betalowerljk$ and $\betaupperljk$.}%
  \label{fig:pooling.problem.network.library.cuts}%
\end{figure}

\noindent
Figure~\ref{fig:pooling.problem.network.library.cuts} shows the variables and sets that, for each pool, output, quality triplet, generate the convex relaxation described in this section.

Linear inequality \eqref{eq:pooling.cut.eq22} is valid for all $(l, j, k) \in L \times J \times K$ where $\betaupperljk > 0$:

\begin{equation}
  \label{eq:pooling.cut.eq22}
  (\etaupperljk - \etalowerljk)\tljk + \etalowerljk(\etaupperljk \slj - \uljk) + \betaupperljk(\uljk - \etalowerljk \slj) \leq \betaupperljk(\pljk  - \etalowerljk)
\end{equation}

\noindent
Linear inequality \eqref{eq:pooling.cut.eq28} is valid for all $(l, j, k) \in L \times J \times K$ where $\betalowerljk < 0$:

\begin{equation}
  \label{eq:pooling.cut.eq28}
  (\etalowerljk - \betalowerljk)(\etaupperljk \slj - \uljk) \leq -\betalowerljk (\etaupperljk - \pljk)
\end{equation}

\noindent
Convex non-linear inequality \eqref{eq:pooling.cut.eq15} is valid  if $\betalowerljk < 0$:

\begin{equation}
  \label{eq:pooling.cut.eq15}
  (\uljk - \betalowerljk \slj)(\uljk - \etalowerljk \slj) \leq -\betalowerljk \slj (\pljk - \etalowerljk)
\end{equation}

\noindent
Convex non-linear inequality \eqref{eq:pooling.cut.eq18} is valid if $\betaupperljk > 0$ and $\etalowerljk < 0$ when $\tljk > 0$:

\begin{equation}
  \label{eq:pooling.cut.eq18}
  (\etaupperljk - \etalowerljk) \tljk + \betaupperljk (\etaupperljk \slj - \uljk) +
  \frac{\etalowerljk \tljk (\uljk - \etalowerljk \slj)}{\tljk + \uljk - \etalowerljk \slj} \leq
  \betaupperljk (\etaupperljk - \pljk)
\end{equation}

Library function \texttt{add\_all\_pooling\_inequalities(block, parent, network)}, which should only be called once, adds the pooling inequalities variables and equations to \texttt{block}, using flow variables from the \texttt{PoolingPQFormulation} block \texttt{parent} and the pooling problem \texttt{network}. This function also adds the Equation~\eqref{eq:pooling.cut.eq22} \& \eqref{eq:pooling.cut.eq28} inequalities. The function adds the following variables to the Pyomo block:

\begin{enumerate}
\item \texttt{z[l, j]}: the total flow $\zlj$ to output $j$ that does not go through pool $l$,
\item \texttt{s[l, j]}: scaled flow $\slj$ from pool $l$ to output $j$,
\item \texttt{u[j, k, l]}: the contribution $\uljk$ of flow through pool $l$ to the quality $k$ at output $j$,
\item \texttt{y[j, k, l]}: the contribution $\tljk$ of the total flow not going through pool $l$ to the quality $k$ at output $j$,
\item \texttt{t[j, k, l]}: the quality $\pljk$ of attribute $k$ going through pool $l$.
\end{enumerate}

The function uses Coramin to add the McCormick envelope of $\rljk = \slj \pljk$ described in Equation~\eqref{eq:pooling.cut.uxp} to the Pyomo block. Finally, the function creates an \texttt{inequalities} \texttt{ConstraintList} that contains all the inequalities from Equation~\eqref{eq:pooling.cut.eq22} \& \eqref{eq:pooling.cut.eq28}. The function iterates over all the combinations of pool, output, quality in the network and checks the conditions of $\betalowerljk$ and $\betaupperljk$ to decide whether or not to add the inequality to the list.
This function adds all possible inequalities. Alternatively, we could add the inequalities as cutting planes, e.g.\ if a relaxation violates the inequalities by more than a parameter $\varepsilon$. Our computational tests on the large scale instances of \citet{Dey2015} and \citet{luedtke2020strong} indicate that the change in performance is not significant.

Equations~\eqref{eq:pooling.cut.eq15} \& \eqref{eq:pooling.cut.eq18} are quadratic and cannot be added directly to the linear relaxation of the pooling problem. The relaxation adds a gradient inequality based on Equation~\eqref{eq:pooling.cut.eq15} for every $(l, j, k) \in L \times J \times K$ where Equation~\eqref{eq:pooling.cut.eq15} is violated more than a parameter $\varepsilon$.
The relaxation adds a gradient inequality based on Equation~\eqref{eq:pooling.cut.eq18} for every $(l, j, k) \in L \times J \times K$ where Equation~\eqref{eq:pooling.cut.eq18} is violated more than a parameter $\varepsilon$ and $\tljk > 0$.

The library function \texttt{add\_valid\_cuts(block, parent, network)} 
checks if any of the non-linear equations, i.e.\ Equations~\eqref{eq:pooling.cut.eq15} \& \eqref{eq:pooling.cut.eq18}, is violated and, if violated, generates a gradient inequality. This function can be called multiple times, e.g.\ in each iteration of a cut loop. 
This function iterates over all the pool, output, quality triplets and adds the cut to the \texttt{cuts} \texttt{ConstraintList} on the Pyomo block if any of the expressions in Equations~\eqref{eq:pooling.cut.eq15} \& \eqref{eq:pooling.cut.eq18} are violated by a parameter $\varepsilon$. Users wishing to generate cuts but not wanting to add the cuts to the model can use the function \texttt{generate\_valid\_cuts(block, parent, network)} that returns an iterator of cuts to be added.

The \texttt{PoolingQPFormulation} block also contains two convenience methods to make the pooling network library more user-friendly:

\begin{enumerate}
\item \texttt{add\_inequalities()}: calls the \texttt{add\_all\_pooling\_inequalities} function with the correct arguments to add the relaxation variables for each pool, output, quality triplet. This method also adds Equation~\eqref{eq:pooling.cut.eq22} \& \eqref{eq:pooling.cut.eq28} if $\betaupperljk > 0$ and $\betalowerljk < 0$ respectively,
\item \texttt{add\_cuts()}: calls the \texttt{add\_valid\_cuts} function to add gradient inequalities based on Equation~\eqref{eq:pooling.cut.eq15} \& \eqref{eq:pooling.cut.eq18}.
\end{enumerate}

Listing~\ref{lst:pooling.network.library.add.cuts} in \ref{app:code} uses these methods to add valid inequalities to the linear relaxation of the pooling problem. This example also shows how to use the \texttt{add\_cuts} method to generate cuts at every iteration of a simple cut loop.

\subsection{Mixed-integer programming restriction}
\label{sect:network.primal.heuristic}

The library implements a MIP restriction of the pooling problem \citep{Dey2015}.
The heuristic ``splits'' each problem pool $l$ in $\tau$ copies, each copy receives a predefined fraction of the total flow received by pool $l$. Each copy of the pool can send flow only to one output node. Equation~\eqref{eq:pooling.mip.restriction} describes the MIP restriction, which introduces auxiliary flow variables $\wiltj$ representing the flow from input $i$ to output $j$ through pool $l$ and auxiliary pool $t$. Each auxiliary pool $t$ receives a fraction $\gammalt$ of the total flow into pool $l$. The restriction also introduces binary variables $\zltj$ that equal $1$ if the auxiliary pool $t$ of pool $l$ sends flow to output $j$, since the heuristic requires that each auxiliary pool sends flow to at most one output, we require that $\sum_{j \in J} \zltj = 1, \forall l \in L, t \in \{1,\dots,\tau\}$.

\begin{subequations}
  \label{eq:pooling.mip.restriction}
  \begin{align}
    & & & \underset{y,v,q}{\text{max}} \sum_{i \in I, j \in J} f_{ij} y_{ij} + \sum_{i \in I, l \in L, j \in J} (f_{il} + f_{lj}) v_{ilj}\\
    & & \begin{array}{r} \text{Flow} \\ \text{Balance} \end{array} &\left[ \vilj = \sum_{t = 1}^{\tau} \wiltj \quad \forall i \in I, l \in L, j \in J
                                                                     \vphantom{\sum_{i \in I}}  \right. \label{eq:pooling.mip.restriction.flow.balance}\\
    & & \begin{array}{r} \text{Fractional Flow} \\ \text{Definition} \end{array} &\left[ \sum_{j \in J} \wiltj = \gammalt \sum_{j \in J} \vilj \quad \forall i \in I, l \in L, t \in \{1,\dots,\tau\} \right. \label{eq:pooling.mip.restriction.fract.flow.def}\\
    & & \begin{array}{r} \text{Bounds} \end{array}&\left[
                                                                    \begin{array}{l}
                                                                      0 \leq \wiltj \leq \clj \zltj \quad \forall i \in I, l \in L, t \in \{1,\dots,\tau\}, j \in J\\
                                                                      \sum_{j \in J} \zltj = 1 \quad \forall l \in L, t \in \{1,\dots,\tau\}\\
                                                                      \zltj \in \{0, 1\} \quad \forall l \in L, t \in \{1,\dots,\tau\}, j \in J
                                                                    \end{array}
    \right.\label{eq:pooling.mip.restriction.extra}\\
    & & \begin{array}{r} \text{PQ-formulation} \end{array}&\left[\text{Equation~\eqref{eq:pooling.pq.simplex} -- Equation~\eqref{eq:pooling.pq.bounds}} \right. \notag
  \end{align}
\end{subequations}

\begin{figure}%
  \centering
  \subfloat[\centering Example pooling problem with 3 inputs, 2 pools, and 2 outputs.]{{\includegraphics[width=0.5\textwidth]{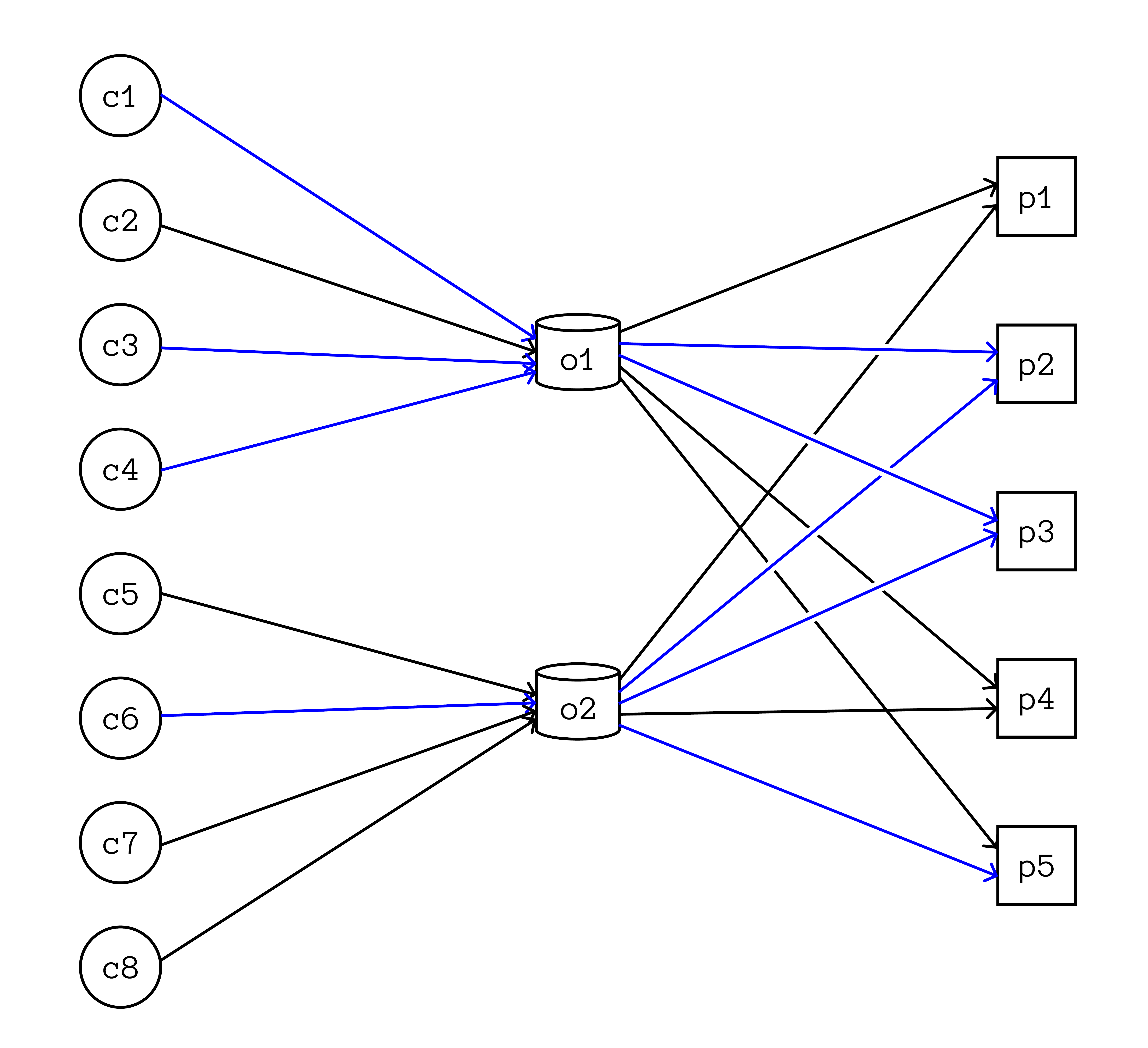} }}%
  \subfloat[\centering MIP restriction of the pooling problem in (a).]{{\includegraphics[width=0.5\textwidth]{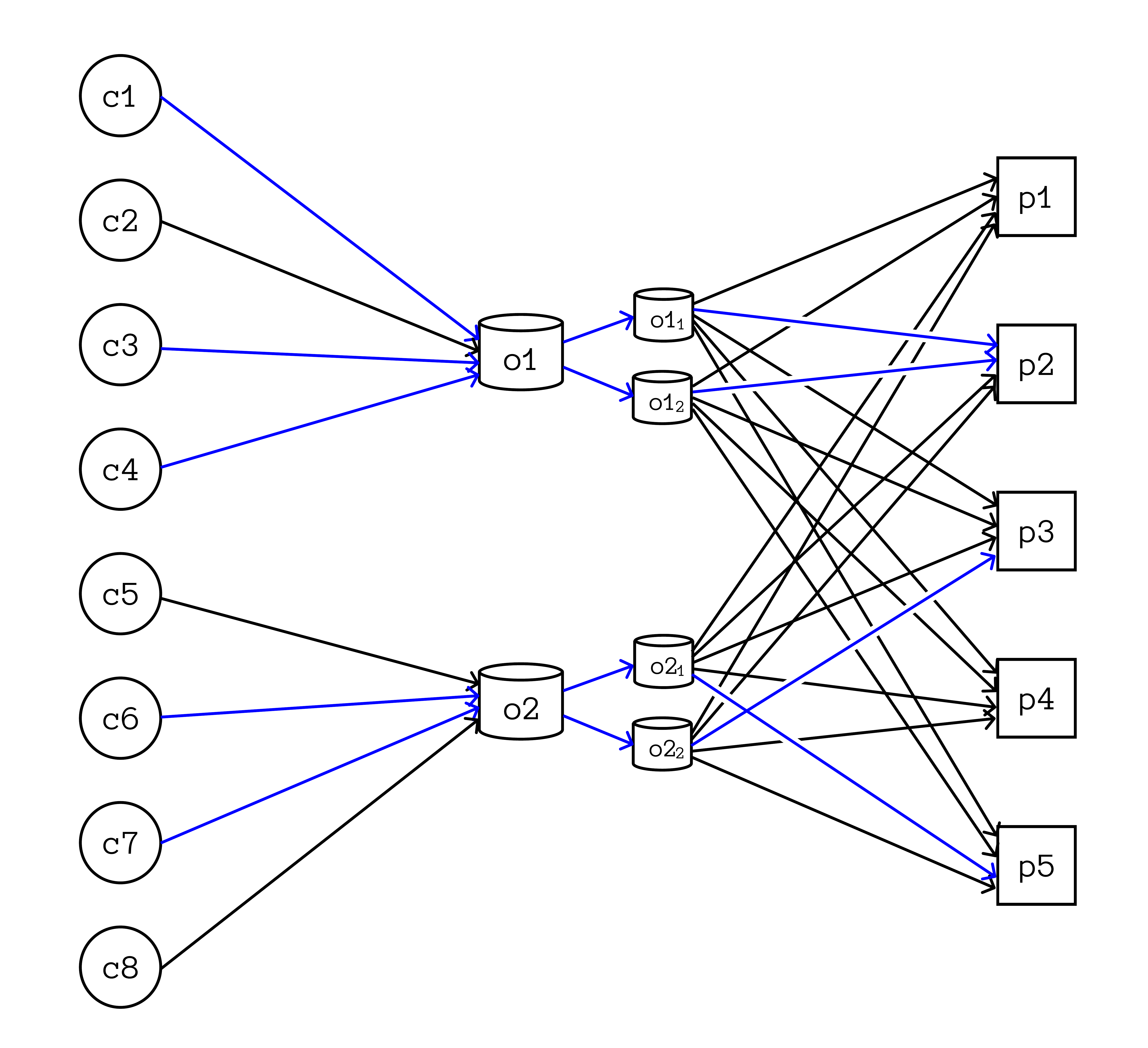}}}%
  \caption{Pooling problem network structure of \texttt{adhya4} (a) together with its MIP restriction (b). The MIP restriction ``splits'' each pool into $\tau$ ($\tau = 2$ in this example) auxiliary pools, each auxiliary pool receives a predetermined fraction of the total flow into its pool and can send flow to only one output. There is no cost associated to the flow between a pool and its auxiliary pools. The blue arrows represent edges with a positive flow in the problem optimal solution and the MIP restriction solution, we can see that in this example the two optimal network configurations are different.}%
  \label{fig:pooling.problem.heuristic}%
\end{figure}

Figure~\ref{fig:pooling.problem.heuristic} diagrams the pooling problem and its MIP restriction.
The heuristic is implemented with two functions: one adds auxiliary variables and constraints to the \texttt{PoolingPQFormulation} block and a second restores the block to its original state. As a convenience to users, the library also implements a Python \textit{context manager} to automatically add and remove the MIP restriction. Listing~\ref{lst:network.add.heuristic} in \ref{app:code} shows how users can build the MIP heuristic.

The context manager \texttt{mip\_heuristic(block, network, tau, weights)} expects the following parameters:

\begin{itemize}
\item \texttt{block}: the \texttt{PoolingPQFormulation} block,
\item \texttt{network}: the pooling problem network,
\item \texttt{tau}: the number $\tau$ of auxiliary pools each pool is ``split'' into,
\item \texttt{weights}: a function to generate the weights $\gammalt$ for each pool $l$ and auxiliary pool $t$. By default, the library uses uniform weights $1/\tau$.
\end{itemize}

The \texttt{mip\_heuristic} context manager modifies the input block by deactivating the block's \texttt{path\_definition} constraints and adding a \texttt{\_pooling\_mip\_heuristic} sub-block with the following variables:

\begin{itemize}
\item \texttt{w[i,l,t,j]}: auxiliary flow variables $\wiltj$,
\item \texttt{zeta[l,t,j]}: binary variables $\zltj$.
\end{itemize}

The sub-block contains the following constraints:

\begin{itemize}
\item \texttt{flow\_balance[i,l,j]}: contains Equation~\eqref{eq:pooling.mip.restriction.flow.balance},
\item \texttt{flow\_balance\_2[i,l,t]}: contains Equation~\eqref{eq:pooling.mip.restriction.fract.flow.def},
\item \texttt{flow\_choice\_limit[i,l,t,j]}: contains the Equation~\eqref{eq:pooling.mip.restriction.extra} $\wiltj$ bounds,
\item \texttt{flow\_choice[l,t]}: limits each auxiliary pool \texttt{t} to send flow to 1 output.
\end{itemize}

Solving the MIP restriction does not produce variable values $\qil$, so the library function \texttt{derive\_fractional\_flow\_variables(block)} derives the $\qil$ values from the other network flows, i.e.\ $\qil = \frac{\vilj}{\ylj} = \frac{\vilj}{\sum_{i \in I}\vilj}$.

\section{\galini{} integration}
\subsection{Initial primal heuristic}
\label{galini.primal.heuristic}

The pooling network library integrates the Section~\ref{sect:network.primal.heuristic} MIP heuristic by providing a \texttt{PoolingPrimalSearchStrategy} class implementing the \galini{} \texttt{InitialPrimalSearchStrategy} interface and registering it with \galini{} to be available at runtime. \galini{} uses this interface to find an initial feasible solution at the root node only.
The primal search strategy detects the Pyomo block that contains the pooling problem and builds the MIP restriction using $\tau = 1$ 
\citep{Dey2015}. 
The search strategy then solves the MIP restriction using the same MIP solver used by the \galini{} branch \& cut algorithm. To limit the time spent finding the feasible solution, the search strategy sets the MIP solver time limit option to $60$ seconds and the optimality gap to $1\%$.

If the MIP heuristic is successful, the search strategy derives the values for the fractional flow variables $\qil$ and then it returns the feasible solution to the original user model to \galini{}.

The primal heuristic is manually activated by setting \galini{} configuration \texttt{intial\_search\_strategy} to \texttt{'pooling'}. \texttt{PoolingPrimalSearchStrategy} does not perform additional checks on the model's additional constraints. If the additional constraints are linear, the primal search strategy will succeed. If the model contains additional non-linear constraints, the primal search strategy will fail since the MIP solver is not able to handle them and Pyomo will throw an exception. We decided not to disable the non-linear constraints since by doing that the new model would not be a restriction of the original model, creating a subtle bug that is difficult to debug for users.

\subsection{Cuts generator}

The pooling network library implements a \galini{} \texttt{CutsGenerator} to generate valid cuts for the pooling problem at each node of the branch \& bound algorithm. The cut generator implements the following interface methods:

\begin{itemize}
\item \texttt{before\_start\_at\_root(problem, linear\_problem)}: this callback is called before the cut loop at the root node. It checks whether \texttt{linear\_problem} contains a \texttt{PoolingPQFormulation} block, if not sets a variable to disable the cut generator. This callback also adds valid inequalities to \texttt{linear\_problem}.
\item \texttt{generate(problem, linear\_problem, lp\_solution, tree, node)}: this method is called at each iteration of the cut loop to generate new cuts to be added to \texttt{linear\_problem}. The function calls the \texttt{generate\_valid\_cuts} function from the pooling network library and then returns a list of cuts that \galini{} adds to \texttt{linear\_problem}.
\end{itemize}

Users can activate this cut generator by adding \texttt{'pooling'} to the list of cuts generator in the \galini{} configuration file.

The cuts generator initializes the linear relaxation of the pooling problem at the root node.
The pooling flow variable bounds change throughout the branch \& bound algorithm, but our preliminary results suggest that these
changes do not significantly affect the \citet{luedtke2020strong} relaxation. Therefore, we do not recompute the relaxation equations in the branch \& bound tree.

\section{Computational Results}

All tests are performed on a machine with an Intel Core i7-6700HQ CPU and 32Gb of RAM. We use \texttt{pooling\_network} version $1.0.0$ \citep{ceccon2021pooling}, \galini{} version $1.0.1$ \citep{ceccon2021galini} using CPLEX version $12.10$ as the MIP solver and Ipopt \citep{wachter2006} version $3.13.2$ to locally solve NLPs. For comparison, we use GUROBI version $9.1$ together with the Pyomo \texttt{gurobi\_direct} solver. We set the optimality relative gap to $10^{-6}$, the optimality absolute gap to $10^{-8}$, and the time limit to $600$ seconds, all other options are left as default. Additionally, we set Gurobi \texttt{NonConvex} option to $2$. We set GUROBI absolute gap by changing the \texttt{MIPGapAbs} option, and we use the \texttt{MIPGap} option to change the relative gap.

The absolute gap is the \textit{absolute difference} between the best possible objective value and the objective value of the best feasible solution.
The relative gap is the \textit{relative difference} $d(a, b)$ between the best possible objective value $a$ and the objective value of the best feasible solution $b$:

\begin{equation}
  d(a, b) = \begin{cases}
    |a - b|/\max\{|a|, |b|\} & \mbox{if } a \neq 0, b \neq 0\\
    \infty & \mbox{if } |a| = \infty \mbox{ or } |b| = \infty\\
    |a - b|/\varepsilon & \mbox{otherwise}, \\
  \end{cases}
  \label{eq:relative_difference}
\end{equation}

\noindent
where $\varepsilon$ is a small quantity. In the instances studied in this section, the objective values are always non-zero and the choice of $\varepsilon$ is irrelevant.

The pooling network library has functions, available under the \texttt{pooling\_} \texttt{network.intances} namespace, to automatically build the network structure of standard pooling problems used in these computational results. The data set summarized in Table~\ref{tbl:pooling.instances} contains classic small-scale problems \citep{tawarmalani2013convexification}, sparse instances generated based on the Haverly instances \citep{luedtke2020strong}, and large-scale dense instances \citep{Dey2015}. 

We test Gurobi with the following configurations:

\begin{itemize}
\item \textbf{Default}: default Gurobi configuration.
\item \textbf{Disable Cuts}: disable cuts generation by setting the \texttt{Cuts} option to $0$.
\item \textbf{Pooling Cuts}: generate \citet{luedtke2020strong} cuts. The test linearizes the optimization problem with the convex hulls of the bilinear terms $\vilj = \qil \ylj$. This configuration solves the linear relaxation of the pooling problem to generate \citet{luedtke2020strong} cuts with at most $20$ cut loop iterations. Finally, the test reintroduces the bilinear terms and solves the problem again with Gurobi. This test considers the time spent solving the linear relaxations and generating cuts as part of the Gurobi solving time.
\item \textbf{Pooling Cuts Oracle}: this configuration generates the same cuts as the \textbf{Pooling Cuts} test, but the time spent generating the cuts and solving the linear relaxations is excluded. This test assumes Gurobi had access to an oracle providing a list of valid cuts.
\item \textbf{Warm Start}: this configuration uses the MIP restriction of the pooling problem to pass a feasible solution to Gurobi. This test considers the time spent solving the MIP restriction as part of the Gurobi solving time.
\item \textbf{Warm Start Oracle}: this configuration uses the same MIP restriction as the previous test, but the time spent solving the MIP restriction is excluded. This test is equivalent to an oracle that provides a good feasible solution to the pooling problem.
\end{itemize}

We test \galini{} with the following configurations:

\begin{itemize}
\item \textbf{Default}: default \galini{} configuration, without pooling specific cuts generator and heuristic active.
\item \textbf{Pooling Cuts Generator}: enable the pooling cuts generator, use the default initial primal heuristic.
\item \textbf{Primal Heuristic}: enable the pooling initial primal heuristic, use the default cuts generators.
\item \textbf{Pooling Cuts Generator and Primal Heuristic}: enable both the pooling cuts generator and the initial primal heuristic.
\end{itemize}

To understand the default configuration, Table~\ref{tbl:large.scale.gurobi.cuts} presents the number of each cut type generated by Gurobi, together with the average number of branch \& bound nodes visited. 
RLT cuts are the only class of cutting planes that is generated consistently for all pooling problems.

We compare configurations with performance profiles \citep{dolan2002benchmarking}, which measure performance of a set of solver configurations $S$ on a test set  $P$. For each problem $p \in P$ and solver $s \in S$, we define a performance measure $t_{ps}$. This performance measure can be the run time needed to solve the problem to global optimality, or the relative gap after a fixed amount of time. For each problem and solver, we define the performance ratio:

\begin{equation}
  \label{eq:large.scale.pp.perf.ratio}
  r_{ps} = \frac{t_{ps}}{\min\{t_{ps} : s \in S\}}.
\end{equation}

The performance profile $\rho_s(\tau)$ is the probability for solver $s \in S$ that a performance ratio $r_{ps}$ is within a factor $\tau \in \mathrm{R}$ of the best possible ratio. 

\begin{equation}
  \label{eq:large.scale.pp.perf.profile}
  \rho_s(\tau) = \frac{1}{\text{card}(P)}\text{card}\{p \in P : r_{ps} \leq \tau\}
\end{equation}

\begin{figure}
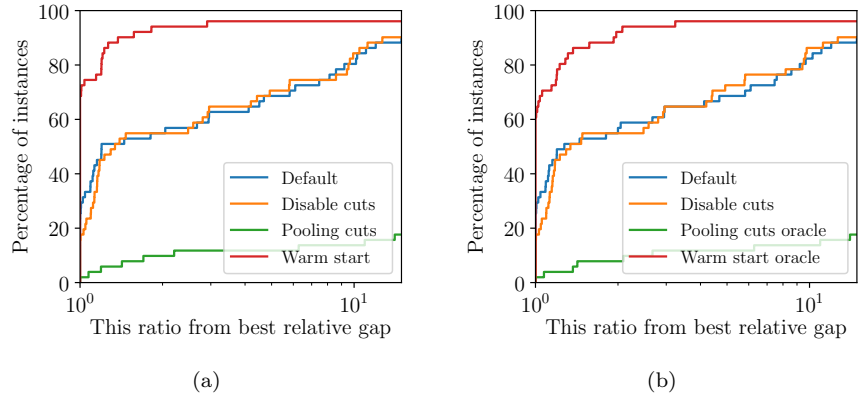

  \centering
  \subfloat[\centering]{
    \resizebox{0.48\textwidth}{!}{%
      \input{pp_gurobi_dense.pgf}
    }
  }
  \subfloat[\centering]{
    \resizebox{0.48\textwidth}{!}{%
      \input{pp_gurobi_dense_oracle.pgf}
    }
  }
  \caption{Performance profile comparing different Gurobi configurations on 50 dense pooling problem instances \citep{Dey2015}. Performance profiles \textbf{a} include the time spent generating the cutting planes and solving the MIP restriction in the total time limit. Performance profiles \textbf{b} use an oracle to generate the cutting planes and find an initial feasible solution.}
  \label{fig:large.scale.pp.gurobi.dense}
\end{figure}

Figure~\ref{fig:large.scale.pp.gurobi.dense} and Table~\ref{tbl:dense.runtime.gurobi} compare  Gurobi configurations on the dense instances \citep{Dey2015}.
Figure~\ref{fig:large.scale.pp.gurobi.dense} shows that Gurobi does not benefit, in these instances, from its own cutting planes: the default performance of Gurobi is similar to the performance when cutting planes are deactivated.
The many \citet{luedtke2020strong} inequalities damage Gurobi's performance, even when adding the cuts as an oracle. 
But the MIP restriction greatly benefits Gurobi, especially as an oracle.
Our experiments confirm the \citet{Dey2015} results: Gurobi is very competitive in solving smaller dense instances (\texttt{randstd11} to \texttt{randstd40}), but for bigger instances a specialized heuristic is beneficial.

Figure~\ref{fig:large.scale.pp.gurobi.sparse} and Table~\ref{tbl:sparse.runtime.gurobi} show the performance of Gurobi configurations on sparse instances  \citep{luedtke2020strong}.
For these instances, the Gurobi default configuration performs best. Adding the pooling cuts before solving the non-convex problem hurts performance because of the overhead of (i) computing cuts and (ii) solving the linear relaxation in the cut loop. The configuration where the pooling cuts are generated as if with an oracle performs similarly to the default Gurobi configuration. The initial primal heuristic is not as beneficial because Gurobi is able to find a good feasible solution by default.
We note that the Gurobi defaults are very strong in solving this class of sparse problems.

\begin{figure}
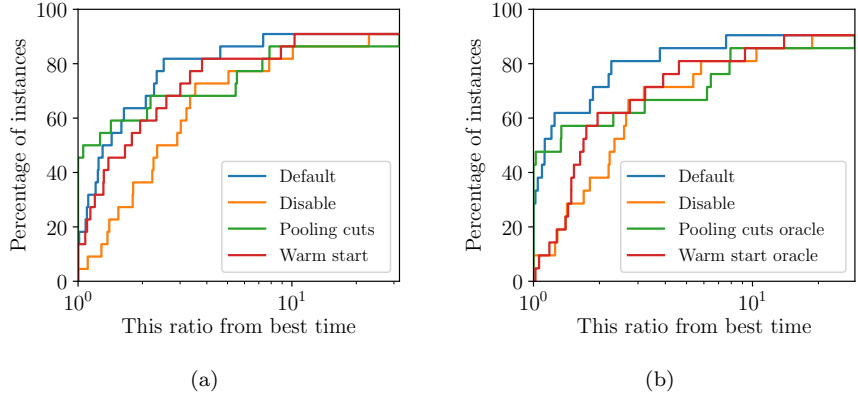

  \centering
  \subfloat[\centering]{
    \resizebox{0.48\textwidth}{!}{%
      \input{pp_gurobi_sparse_big.pgf}
    }
  }
  \subfloat[\centering]{
    \resizebox{0.48\textwidth}{!}{%
      \input{pp_gurobi_sparse_big_oracle.pgf}
    }
  }
  \caption{Performance profile comparing different Gurobi configurations on 30 difficult sparse pooling problem instances~\citep{luedtke2020strong}. We consider an instance difficult if the Gurobi default configuration takes longer than $10$ seconds to solve it to optimality. Performance profiles \textbf{a} include the time spent generating the cutting planes and solving the MIP restriction in the total time limit, while performance profiles \textbf{b} use an oracle to generate the cutting planes and find an initial feasible solution.}
  \label{fig:large.scale.pp.gurobi.sparse.big}
\end{figure}

To consider whether the pooling cuts help solve more difficult sparse problems, we consider the $30$ sparse instances where the Gurobi default configuration requires $>10$ seconds to solve to global optimality. We plot the performance profile on these instances in Figure~\ref{fig:large.scale.pp.gurobi.sparse.big}: for these instances the pooling cuts are beneficial and provide a performance boost over the default Gurobi.
As observed by \citet{luedtke2020strong}, the pooling cuts are not as helpful on denser instances.
In our results the benefit provided by the pooling cuts is not as strong as in~\citet{luedtke2020strong} because Gurobi is more performant than the ones used in the original paper.

\begin{figure}
  \centering
  \subfloat[\centering]{
    \resizebox{0.48\textwidth}{!}{%
      \input{pp_gurobi_sparse.pgf}
    }
  }
  \subfloat[\centering]{
    \resizebox{0.48\textwidth}{!}{%
      \input{pp_gurobi_sparse_oracle.pgf}
    }
  }
  \caption{Performance profile comparing different Gurobi configurations on 360 sparse pooling problem instances \citep{luedtke2020strong}. Performance profiles \textbf{a} include the time spent generating the cutting planes and solving the MIP restriction in the total time limit. Performance profiles \textbf{b} use an oracle to generate the cutting planes and find an initial feasible solution.}
  \label{fig:large.scale.pp.gurobi.sparse}
\end{figure}

\begin{figure}
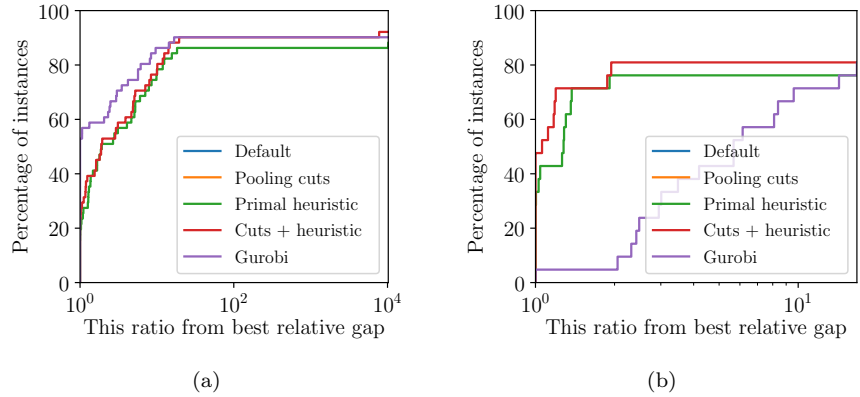

  \centering
  \subfloat[\centering]{
    \resizebox{0.48\textwidth}{!}{%
      \input{pp_galini_dense.pgf}
    }
  }
  \subfloat[\centering]{
    \resizebox{0.48\textwidth}{!}{%
      \input{pp_galini_dense_big.pgf}
    }
  }
  \caption{Performance profile of different \galini{} configurations compared to Gurobi on \textbf{a)} 50 large scale dense instances~\citep{Dey2015}, and \textbf{b)} the 20 largest dense instances. The default and pooling cuts configurations did not find a feasible solution and so they are not visible in the plots.}
  \label{fig:large.scale.pp.galini.dense}
\end{figure}

Figure~\ref{fig:large.scale.pp.galini.dense} and Table~\ref{tbl:dense.runtime.galini} compare the different \galini{} configurations and the default Gurobi configuration on 50 dense pooling instances. Activating the initial primal heuristic improves \galini{} performance and makes it comparable to Gurobi. When we compare \galini{} with Gurobi only on the 20 largest dense instances, we can see that \galini{} with the initial primal heuristic enabled has a very strong performance.
Table~\ref{tbl:large.scale.result.summary} compare the relative gaps between \galini{} (with pooling cuts and initial primal heuristic enabled) and Gurobi on dense instances: the \galini{} relative gaps are stable as the problem sizes increases, while Gurobi relative gap degrades with increasing problems size.

\begin{figure}
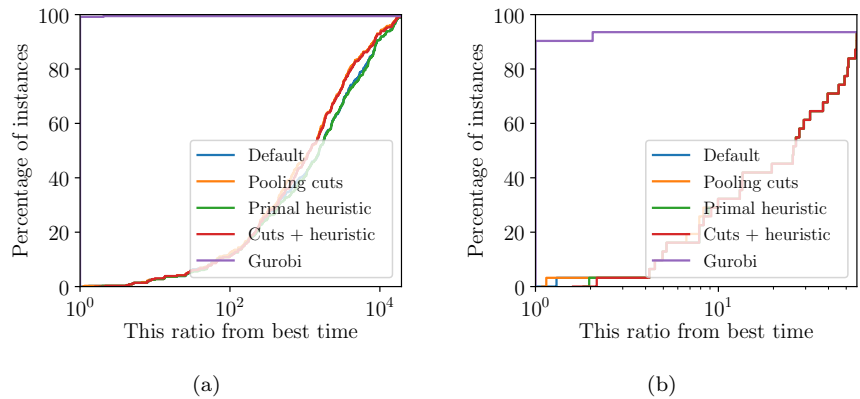

  \centering
  \subfloat[\centering]{
    \resizebox{0.48\textwidth}{!}{%
      \input{pp_galini_sparse.pgf}
    }
  }
  \subfloat[\centering]{
    \resizebox{0.48\textwidth}{!}{%
      \input{pp_galini_sparse_big.pgf}
    }
  }
  \caption{Performance profile of different \galini{} configurations compared to Gurobi on \textbf{a)} 360 large scale sparse instances~\citep{luedtke2020strong}, and \textbf{b)} 30 difficult sparse instances.  The default and pooling cuts configurations did not find a feasible solution and so they are not visible in the plots.}
  \label{fig:large.scale.pp.galini.sparse}
\end{figure}

Figure~\ref{fig:large.scale.pp.galini.sparse} and Table~\ref{tbl:sparse.runtime.galini} compare different \galini{} configurations and the default Gurobi configuration on 360 sparse pooling instances. In this case, the \galini{} performance is not competitive with Gurobi. 
The GALINI branch \& bound algorithm is not as fast as Gurobi's in visiting the hundreds or thousands of nodes required to solve these problems.

To understand the difference in relative gaps between \galini{} and Gurobi in the dense instances, Table~\ref{tbl:large.scale.bounds.difference} shows the relative difference between the lower bounds and upper bounds computed by the two solvers.
The objective lower bounds computed by the two solvers are fairly similar, while the upper bound computed by \galini{} for the larger problems is on average $28\%$ better than the upper bound computed by Gurobi. Therefore, \galini{} performs better on large scale pooling problems because it's able to exploit the pooling problem structure to compute a feasible solution using a pooling-specific heuristic.

\section{Conclusion}

This manuscript presented the pooling network library, a Python library to describe pooling problems and generate Pyomo models using the PQ-formulation. The pooling network library provides functions to add inequalities and cuts to build a convex relaxation of the non-convex problem, which helps in providing a stronger lower bound on the problem solution. The library also includes functions to build a MIP restriction of the pooling problem and find a feasible solution using a MIP solver.

We then showed that, by integrating the cuts generator and heuristic from Section~\ref{chapt:network.structure}, we can solve large scale dense pooling problems to a smaller gap than the best available commercial solver. We compare the lower and upper bounds produced by the two solvers and see that the \galini{} performance derives from a better primal heuristic that takes advantage of the pooling problem structure.
Our experiments also show that, on sparse pooling instances, the MIP heuristic is not as effective since Gurobi already has strong performance on this class of problems. We show that, on dense pooling instances, solving the MIP restriction of the pooling problem can improve Gurobi's performance.

\section{Acknowledgements}
This work was funded by an Engineering \& Physical Sciences Research Council Research Fellowship to RM [GrantNumber EP/P016871/1].

\bibliography{biblio}


\begin{landscape}
\begin{table}
\centering
  \caption{Pooling instances summary. The table contains each class of problems, together with a summary on the network size and the number of variables in the PQ-formulation of the problem. Sources: A) \citet{luedtke2020strong}, B) \citet{Dey2015}. All instances are available under the \texttt{pooling\_network.instances} namespace.}
  \label{tbl:pooling.instances}
  \begin{tabular}{lcccccccccc}
    \toprule
    \multicolumn{7}{c}{} & \multicolumn{3}{c}{\underline{Avg. no. vars.}}\\
    Problem & Source & Instances & $|I|$ & $|L|$ & $|J|$ & $|K|$ & $|A|$ & $\qil$ & $\ylj$ & $\yij$ \\
    \midrule
    \texttt{haverly\_10\_ae\_10\_*}& \multirow{17}{*}{A} & 20 & 30 & 10 & 20 & 1/2 & 70 & 22 & 22 & 26\\
    \texttt{haverly\_10\_ae\_20\_*}& & 20 & 30 & 10 & 20 & 1/2 & 80 & 25 & 24 & 31\\
    \texttt{haverly\_10\_ae\_30\_*}& & 20 & 30 & 10 & 20 & 1/2 & 90 & 27 & 25 & 38\\
    \texttt{haverly\_10\_ae\_40\_*}& & 20 & 30 & 10 & 20 & 1/2 & 100 & 32 & 26 & 42\\
    \texttt{haverly\_10\_ae\_50\_*}& & 20 & 30 & 10 & 20 & 1/2 & 110 & 32 & 29 & 49\\
    \texttt{haverly\_10\_ae\_60\_*}& & 20 & 30 & 10 & 20 & 1/2 & 120 & 35 & 29 & 56\\
    \texttt{haverly\_15\_ae\_15\_*}& & 20 & 45 & 15 & 30 & 1/2 & 105 & 33 & 33 & 38\\
    \texttt{haverly\_15\_ae\_30\_*}& & 20 & 45 & 15 & 30 & 1/2 & 120 & 39 & 35 & 47\\
    \texttt{haverly\_15\_ae\_45\_*}& & 20 & 45 & 15 & 30 & 1/2 & 130 & 42 & 38 & 55\\
    \texttt{haverly\_15\_ae\_60\_*}& & 20 & 45 & 15 & 30 & 1/2 & 150 & 45 & 40 & 65\\
    \texttt{haverly\_15\_ae\_75\_*}& & 20 & 45 & 15 & 30 & 1/2 & 165 & 50 & 45 & 70\\
    \texttt{haverly\_15\_ae\_90\_*}& & 20 & 45 & 15 & 30 & 1/2 & 180 & 54 & 46 & 80\\
    \texttt{haverly\_20\_ae\_20\_*}& & 20 & 60 & 20 & 40 & 1/2 & 140 & 46 & 44 & 50\\
    \texttt{haverly\_20\_ae\_40\_*}& & 20 & 60 & 20 & 40 & 1/2 & 160 & 52 & 47 & 62\\
    \texttt{haverly\_20\_ae\_60\_*}& & 20 & 60 & 20 & 40 & 1/2 & 180 & 56 & 50 & 74\\
    \texttt{haverly\_20\_ae\_80\_*}& & 20 & 60 & 20 & 40 & 1/2 & 200 & 61 & 55 & 84\\
    \texttt{haverly\_20\_ae\_100\_*}& & 20 & 60 & 20 & 40 & 1/2 & 220 & 68 & 56 & 96\\
    \texttt{haverly\_20\_ae\_120\_*}& & 20 & 60 & 20 & 40 & 1/2 & 240 & 72 & 61 & 108\\
    \hline
    \texttt{randstd11-20}& \multirow{5}{*}{B} & 10 & 25 & 18 & 25 & 8 & 408 & 202 & 183 & 23\\
    \texttt{randstd21-30}& & 10 & 25 & 22 & 30 & 10 & 542 & 240 & 273 & 29\\
    \texttt{randstd31-40}& & 10 & 30 & 22 & 35 & 10 & 694 & 289 & 230 & 40\\
    \texttt{randstd41-50}& & 10 & 40 & 30 & 45 & 10 & 1136 & 500 & 566 & 70\\
    \texttt{randstd51-60}& & 10 & 40 & 30 & 50 & 14 & 1209 & 501 & 631 & 77\\
    \bottomrule
  \end{tabular}
\end{table}
\end{landscape}

\begin{landscape}
\begin{table}
\centering
  \caption{Average number of cuts generated and nodes explored by Gurobi with the default configuration. Empty cells mean that no cuts were generated, while cell with value $0$ mean that some cuts were generated but the average is less than $1$.}
  \label{tbl:large.scale.gurobi.cuts}
  \begin{tabular}{lcccccccc}
    \toprule
    Problem & Gomory & Implied bound & Proj. implied bound & MIR & Flow cover & RLT & Relax-and-lift & No. nodes\\
    \midrule
    \texttt{haverly\_10\_ae\_10\_*}& & & & & & 5 & & 168 \\
    \texttt{haverly\_10\_ae\_20\_*}& & & & & & 8 & & 609 \\
    \texttt{haverly\_10\_ae\_30\_*}& & & & & & 6 & & 636 \\
    \texttt{haverly\_10\_ae\_40\_*}& & & & & & 5 & & 1196 \\
    \texttt{haverly\_10\_ae\_50\_*}& & & & & & 5 & & 2822 \\
    \texttt{haverly\_10\_ae\_60\_*}& & & & & & 7 & & 2092 \\
    \texttt{haverly\_15\_ae\_15\_*}& & & & & & 11 & & 885 \\
    \texttt{haverly\_15\_ae\_30\_*}& & & & & & 10 & & 2694 \\
    \texttt{haverly\_15\_ae\_45\_*}& & & & & & 8 & & 14603 \\
    \texttt{haverly\_15\_ae\_60\_*}& & & & & & 5 & & 49238 \\
    \texttt{haverly\_15\_ae\_75\_*}& & & & & & 7 & & 9366 \\
    \texttt{haverly\_15\_ae\_90\_*}& & & & & & 4 & & 91873 \\
    \texttt{haverly\_20\_ae\_20\_*}& & & & & & 8 & & 2378 \\
    \texttt{haverly\_20\_ae\_40\_*}& & & & & & 9 & & 7704 \\
    \texttt{haverly\_20\_ae\_60\_*}& & & & & & 5 & & 50797 \\
    \texttt{haverly\_20\_ae\_80\_*}& & & & & & 6 & & 45601 \\
    \texttt{haverly\_20\_ae\_100\_*}& & & & & & 4 & & 34336 \\
    \texttt{haverly\_20\_ae\_120\_*}& & & & & & 3 & & 52359 \\
    \hline
    \texttt{randstd11-20}& 1 & 1 & 0 & 4 & 7 & 26 & 0 & 62362 \\
    \texttt{randstd21-30}& 0 & 0 & 0 & 0 & 1 & 15 & 0 & 37326 \\
    \texttt{randstd31-40}& 2 & 0 & 0 & 0 & 2 & 18 & 0 & 15291 \\
    \texttt{randstd41-50}& 0 & 2 &   & 1 & 1 & 53 & 0 & 1556 \\
    \texttt{randstd51-60}& 0 & 0 &   & 0 & 1 & 42 & 0 & 2400 \\
    \bottomrule
  \end{tabular}
\end{table}
\end{landscape}

\begin{landscape}
\begin{table}
\centering
  \caption{This table compares the shifted geometric mean (shift = 0.1) of the relative gap after 600 seconds of different Gurobi configurations. The gap column shows the relative gap in percentage, the count column the number of instances where the solver is able to find a feasible solution. Instances are from~\citep{Dey2015} and are grouped by size.}
  \label{tbl:dense.runtime.gurobi}
  \begin{tabular}{l cc cc cc cc cc cc}
    \toprule
    \texttt{randstd} & \multicolumn{2}{c}{Default} & \multicolumn{2}{c}{Disable cuts} & \multicolumn{2}{c}{Pooling cuts} & \multicolumn{2}{c}{Warm start} & \multicolumn{2}{c}{Pooling cuts oracle} & \multicolumn{2}{c}{Warm start oracle}\\
    Instance & Gap (\%) & Count & Gap (\%) & Count & Gap (\%) & Count & Gap (\%) & Count & Gap (\%) & Count & Gap (\%) & Count\\
    \midrule
    \texttt{11-20}&0.87	&9&	1.01&	9&	1.70&	5&	2.18&	10&	1.84&	5&	2.38&	10\\
    \texttt{21-30}&1.38	&10&	1.50&	10&	1.35&	3&	1.39&	10&	1.32&	3&	1.39&	10\\
    \texttt{31-40}&1.81&	10&	1.55&	10&	5.64&	2&	1.47&	10&	5.64&	2&	1.52&	10\\
    \texttt{41-50}&38.44&	10	&36.89&	10&&		0&	10.77&	10&&		0	&10.83&	10\\
    \texttt{51-60}&13.44&	8&	11.84&	8&&		0&	2.08&	10&&		0&	2.08&	10\\
    \bottomrule
  \end{tabular}
\end{table}
\end{landscape}

\begin{landscape}
\begin{table}
\centering
  \caption{This table shows the shifted geometric mean (shift = 2) of the run time (in seconds) of different Gurobi configurations. Instances are from ~\citep{luedtke2020strong} and are grouped by size.}
  \label{tbl:sparse.runtime.gurobi}
  \begin{tabular}{l cc cc cc cc cc cc}
    \toprule
    \texttt{haverly} & \multicolumn{2}{c}{Default} & \multicolumn{2}{c}{Disable cuts} & \multicolumn{2}{c}{Pooling cuts} & \multicolumn{2}{c}{Warm start} & \multicolumn{2}{c}{Pooling cuts oracle} & \multicolumn{2}{c}{Warm start oracle}\\
    Instance & Time (s) & Count & Time (s) & Count & Time (s) & Count & Time (s) & Count & Time (s) & Count & Time (s) & Count\\
    \midrule
    \texttt{10\_ae\_10\_*}& 0.05 & 20 &	0.06 &	20 &	1.04 &	20 &	0.15 &	20 &	0.06 &	20 &	0.06 &	20\\
    \texttt{10\_ae\_20\_*}&0.12	&20 &	0.18 &	20 &	1.37 &	20 &	0.21	& 20 &	0.22	&20 &	0.12	&20\\
    \texttt{10\_ae\_30\_*}& 0.14	&20&	0.22	& 20 &	1.12	& 20 &	0.26&	20&	0.14	&20&	0.16	&20\\
    \texttt{10\_ae\_40\_*}& 0.21	&20	& 0.24&	20&	1.96&	20&	0.37	&20&	0.28	&20 &	0.24&	20\\
    \texttt{10\_ae\_50\_*}& 0.31	&20&	0.32	&20 &	1.49	&20 &	0.50	&20	& 0.54&	20	& 0.42 &	20\\
    \texttt{10\_ae\_60\_*}& 0.35	&20&	0.43&	20&	1.83	& 20&	0.50&	20&	0.43&	20&	0.47&	20\\
    \texttt{15\_ae\_15\_*}& 0.16&	20&	0.19&	20&	1.85&	20&	0.26&	20&	0.18&	20	& 0.12&	20\\
    \texttt{15\_ae\_30\_*}& 0.40&	20&	0.88&	20	& 2.15&	20&	0.54	&20& 	0.41&	20&	0.48&	20\\
    \texttt{15\_ae\_45\_*}& 1.00	&20&	1.45	&20&	2.82&	20&	0.99&	20&	0.91&	20&	0.91&	20\\
    \texttt{15\_ae\_60\_*}&2.20&	20	&2.85&	20&	4.48&	20&	2.55&	20&	2.41&	20&	2.79&	20\\
    \texttt{15\_ae\_75\_*}& 1.39&	20&	1.89&	20&	4.51&	20&	1.77&	20&	2.20&	20	& 1.85&	20\\
    \texttt{15\_ae\_90\_*}& 5.58&	20&	6.82&	20&	9.82&	20&	6.39&	20&	7.08&	20&	7.09&	20\\
    \texttt{20\_ae\_20\_*}& 0.34&	20&	1.32&	20&	3.14&	20&	0.92&	20&	0.52&	20&	0.78&	20\\
    \texttt{20\_ae\_40\_*}&1.13&	20&	3.11&	20	& 3.94&	20&	1.46&	20&	1.00&	20&	1.48&	20 \\
    \texttt{20\_ae\_60\_*}&4.79&	20	& 6.94&	20 &	9.40&	20&	5.27&	20	&6.27&	20&	5.97&	20 \\
    \texttt{20\_ae\_80\_*}& 4.12&	20&	6.41&	20	& 9.40&	20&	4.52&	20&	7.49&	20&	5.04&	20\\
    \texttt{20\_ae\_100\_*}& 4.05&	20&	4.62&	20&	8.94&	20&	4.73&	20&	7.00&	20&	5.34&	20\\
    \texttt{20\_ae\_120\_*}& 5.62&	20	&8.83&	20&	8.70&	20&	6.38&	20&	7.00&	20	&6.94	&20\\
    \bottomrule
  \end{tabular}
\end{table}
\end{landscape}

\begin{landscape}
\begin{table}
\centering
  \caption{This table compares the shifted geometric mean (shift = 0.1) of the relative gap after 600 seconds of different GALINI configurations and Gurobi. The gap column shows the relative gap in percentage, the count column the number of instances where the solver is able to find a feasible solution. Instances are from~\citep{Dey2015} and are grouped by size.}
  \label{tbl:dense.runtime.galini}
  \begin{tabular}{l cc cc cc cc cc}
    \toprule
    \texttt{randstd} & \multicolumn{2}{c}{Default} & \multicolumn{2}{c}{Pooling cuts} & \multicolumn{2}{c}{Primal heuristic} & \multicolumn{2}{c}{Cuts + heuristic} & \multicolumn{2}{c}{Gurobi}\\
    Instance & Gap (\%) & Count & Gap (\%) & Count & Gap (\%) & Count & Gap (\%) & Count & Gap (\%) & Count\\
    \midrule
    \texttt{11-20}& &	0&&		0&	4.40&	9&	4.70&	10&	0.84&	9\\
    \texttt{21-30}&&	0	&&	0&	3.95&	10&	3.93&	10&	1.31&	10\\
    \texttt{31-40}&&	0&&		0&	3.21&	10&	3.13&	10&	1.68&	10\\
    \texttt{41-50}&&	0&&		0&	8.57&	9&	8.83&	9&	36.65&	10\\
    \texttt{51-60}&&	0&&		0&	3.31&	8&	3.36&	9&	12.71&	8\\
    \bottomrule
  \end{tabular}
\end{table}
\end{landscape}

\begin{table}
\centering
  \caption{Breakdown of the differences between \galini{} and Gurobi by problem size. \galini{} results are stable across problem sizes, while Gurobi relative gap degrades as the problem size increases. The last three columns show the number of times \galini{} has the best lower bound, upper bound, and relative gap out of the 10 instances in each group.}
  \label{tbl:large.scale.result.summary}
  {\small
  \begin{tabular}{l|ccc|ccc|ccc}
    \toprule
    \texttt{randstd} & \multicolumn{3}{c}{\galini{} Rel. Gap \%} & \multicolumn{3}{c}{Gurobi Rel. Gap \%} & \multicolumn{3}{c}{\galini{} no. wins}\\
    Instance & min & avg & max & min & avg & max & LB & UB & Rel. Gap\\
    \midrule
    \texttt{11-20}& 1.47 & 5.64 & 16.12 & 0.19 & 0.87 & 2.98 & 2 & 1 & 1\\
    \texttt{21-30}& 0.78 & 4.37 & 8.24 & 0.0 & 1.39 & 4.38 & 0 & 0 & 0\\
    \texttt{31-40}& 1.78 & 3.31 & 5.77 & 0.12 & 1.82 & 6.69 & 3 & 3 & 2\\
    \texttt{41-50}& 1.39 & 10.47 & 20.63 & 18.76 & 38.88 & 66.05 & 1 & 9 & 9\\
    \texttt{51-60}& 1.22 & 3.60 & 6.54 & 1.03 & 13.50 & 19.36 & 1 & 9 & 9\\
    \bottomrule
  \end{tabular}
  }
\end{table}

\begin{landscape}
\begin{table}
\centering
  \caption{This table shows the shifted geometric mean (shift = 2) of the run time (in seconds) of different GALINI configurations and Gurobi. Instances are from ~\citep{luedtke2020strong} and are grouped by size.}
  \label{tbl:sparse.runtime.galini}
  \begin{tabular}{l cc cc cc cc cc}
    \toprule
    \texttt{haverly} & \multicolumn{2}{c}{Default} & \multicolumn{2}{c}{Pooling cuts} & \multicolumn{2}{c}{Primal heuristic} & \multicolumn{2}{c}{Cuts + heuristic} & \multicolumn{2}{c}{Gurobi}\\
    Instance & Time (s) & Count & Time (s) & Count & Time (s) & Count & Time (s) & Count & Time (s) & Count\\
    \midrule
    \texttt{10\_ae\_10\_*}&365.08&	20&	229.59&	20&	368.38&	20&	243.90&	20&	0.05&	20\\
    \texttt{10\_ae\_20\_*}& 430.20&	20&	286.16&	20&	458.51&	20&	298.89&	20&	0.12&	20\\
    \texttt{10\_ae\_30\_*}&431.19&	20&	263.31&	20&	449.12&	20&	265.55&	20&	0.14&	20\\
    \texttt{10\_ae\_40\_*}&503.12&	20&	261.69&	20&	519.11&	20&	260.91&	20&	0.21&	20\\
    \texttt{10\_ae\_50\_*}&383.95&	20&	318.31&	20&	396.01&	20&	330.40&	20&	0.31&	20\\
    \texttt{10\_ae\_60\_*}&565.30&	20&	529.80&	20&	574.95&	20&	542.18&	20&	0.35&	20\\
    \texttt{15\_ae\_15\_*}&443.14&	19&	377.27&	20&	498.70&	20&	428.15&	20&	0.16&	20\\
    \texttt{15\_ae\_30\_*}&490.03&	20&	331.77&	20&	505.00&	20&	338.67&	20&	0.40&	20\\
    \texttt{15\_ae\_45\_*}&386.64&	20&	328.40&	20&	408.09&	20&	340.00&	20&	1.00&	20\\
    \texttt{15\_ae\_60\_*}&543.69&	20&	513.30&	20&	586.86&	20&	542.55&	20&	2.20&	20\\
    \texttt{15\_ae\_75\_*}&508.98&	20&	505.21&	20&	526.85&	20&	535.45&	20&	1.39&	20\\
    \texttt{15\_ae\_90\_*}&599.87&	20&	599.99&	20&	600.06&	20&	599.70&	20&	5.58&	20\\
    \texttt{20\_ae\_20\_*}&420.50&	20&	380.41&	20&	468.82&	20&	394.99&	19&	0.34&	20\\
    \texttt{20\_ae\_40\_*}&453.72&	20&	435.21&	20&	485.98&	20&	477.57&	20&	1.13&	20\\
    \texttt{20\_ae\_60\_*}&467.58&	20&	459.88&	20&	509.36&	19&	502.02&	20&	4.79&	20\\
    \texttt{20\_ae\_80\_*}&483.04&	20&	464.25&	20&	494.04&	20&	489.85&	20&	4.12&	20\\
    \texttt{20\_ae\_100\_*}&575.08&	20&	575.59&	20&	586.58&	20&	587.44&	20&	4.05&	20\\
    \texttt{20\_ae\_120\_*}&569.00&	20&	565.16&	20&	582.06&	20&	573.33&	20&	5.62&	20\\
    \bottomrule
  \end{tabular}
\end{table}
\end{landscape}

\begin{table}
\centering
  \caption{Breakdown of the differences between \galini{} and Gurobi lower and upper bounds by problem size. The table contains the relative difference between each solver own lower bound (upper bound) and the best lower bound, averaged between the 10 instances in each group. An average relative difference of $0.00\%$ means the solver always has the best lower or upper bound in that problem size group.}
  \label{tbl:large.scale.bounds.difference}
  {\small
  \begin{tabular}{l|cc|cc}
    \toprule
    \texttt{randstd} & \multicolumn{2}{c}{Avg. Lower Bounds Rel. Diff. \%} & \multicolumn{2}{c}{Avg. Upper Bounds Rel. Diff. \%}\\
    Instance & \galini{} & Gurobi & \galini{} & Gurobi\\
    \midrule
    \texttt{11-20}& 0.27 & 0.10 & 2.84 & 0.00\\
    \texttt{21-30}& 0.33 & 0.00 & 2.48 & 0.00\\
    \texttt{31-40}& 0.09 & 0.00 & 1.47 & 0.02\\
    \texttt{41-50}& 0.13 & 0.00 & 0.00 & 28.32\\
    \texttt{51-60}& 0.12 & 0.00 & 0.00 & 9.00\\
    \bottomrule
  \end{tabular}
  }
\end{table}


\pagebreak
\appendix

\section{Notation}\label{app:notation}

Tables~\ref{table:pooling.problem.notation.sets} and \ref{table:pooling.problem.notation.variables} summarize the notation in this manuscript.

\begin{table}[h]
  \begin{center}
    \begin{tabularx}{\textwidth}{@{} c c X @{}}
      \toprule
      Type & Name & Description\\
      \midrule
      Sets & $N$ & The set of nodes in the network, $N = I \cup L \cup J$\\
      & $I$ & The set of inputs in the network\\
      & $L$ & The set of pools in the network\\
      & $J$ & The set of outputs in the network\\
      & $A$ & The set of edges between nodes in the network\\
      & $I_l$ & The set of inputs connected to pool $l \in L$\\
      & $I_j$ & The set of inputs connected to output $j \in J$\\
      & $K$ & The set of material qualities\\
      \midrule
      Indices & $i \in I$ & Input streams (raw materials or feed stocks)\\
      & $l \in L$ & Pools (blending facilities)\\
      & $j \in J$ & Output streams (end products)\\
      & $k \in K$ & Attributes (qualities monitored)\\
      & $t \in \{1, 2, ..., \tau\}$ & Auxiliary pools\\
      \bottomrule
    \end{tabularx}
  \end{center}
  \caption{Notation for pooling problems used in this manuscript: sets and indices.}
  \label{table:pooling.problem.notation.sets}
\end{table}

\begin{table}[h]
  \begin{center}
    \begin{tabularx}{\textwidth}{@{} c c X @{}}
      \toprule
      Type & Name & Description\\
      \midrule
      Variables & $\qil$ & Fractional flow from input $i$ to pool $l$, as fraction of total incoming flow into pool $l$\\
      & $\ylj$ & Flow from input $i$ to pool $l$\\
      & $\yij$ & Bypass flow from input $i$ to output $j$\\
      & $(of_j)$ & Outflow of product $j$\\
      & $\plk$ & Level of quality attribute $k$ in pool $l$\\
      & $\zlj$ & Total flow to output $j$ that does not pass through pool $l$\\
      & $\uljk$ & Contribution of flow through pool $l$  to quality $k$ at output $j$\\
      & $\tljk$ & Contribution of the total flow not going through pool $l$ to the quality $k$ at output $j$\\
      & $\pljk$ & Quality of attribute $k$ at pool $l$\\
      & $\slj$ & Scaled flow variable of flow between pool $l$ and output $j$\\
      & $\rljk$ & Auxiliary variable, $\rljk = \slj \pljk$\\
      & $\wiltj$ & Flow from input $i$ to output $j$ through pool $l$ and auxiliary pool $t$\\
      & $\zltj$ & Binary variable that equals $1$ if the auxiliary pool $t$ sends flow to output $j$\\      
      \midrule
      Parameters & $\xil$ & Flow from input $i$ to pool $l$\\
      & $\ci$ & Unit cost of raw material feed stock $i$\\
      & $\rj$ & Unit revenue of product $j$\\
      & $\Ail$ -- $\Aiu$ & Availability bounds of input $i$\\
      & $\Sl$ & Capacity of pool $l$\\
      & $\Djl$ -- $\Dju$ & Demand bounds for product $j$\\
      & $\Cik$ & Level of quality $k$ in raw material feed stock $i$\\
      & $\Pjkl$ -- $\Pjku$ & Acceptable composition range of quality $k$ in product $j$\\
      & $\etalowerljk$ -- $\etaupperljk$ & Lower and upper bounds on the excess of quality $k$ at output $j$ over the inputs $i$ that are connected to pool $l$\\
      & $\betalowerljk$ -- $\betaupperljk$ & Lower and upper bounds on the excess of quality $k$ at output $j$ over the inputs $i$ that are \textit{not} connected to pool $l$\\
      & $\gammalt$ & Fraction of the total flow into pool $l$ that the auxiliary pool $t$ receives\\
    \bottomrule
    \end{tabularx}
  \end{center}
  \caption{Notation for pooling problems used in this manuscript: variables and parameters.}
  \label{table:pooling.problem.notation.variables}
\end{table}

\section{Pyomo smart block example}\label{app:pyomo_smart_block}

Listing~\ref{lst:pyomo.example.2} in \ref{app:code} shows how to define a new \texttt{McCormickEnvelope} block that automatically generates the four inequalities of the McCormick envelope of $w = x y$ \citep{McCormick1976}. To define our block data, we define a new class that subclasses Pyomo \texttt{\_BlockData}, and call \texttt{\_BlockData.\_\_init\_\_(component)} method from our class initialization function. To define our block type, we use the \texttt{@declare\_custom\_block(name)} decorator on our class. A decorator is a special type of Python function that ``wraps'' another function or, like in this case, class. The name we pass to the \texttt{declare\_custom\_block} decorator is the name of our block type. The name of the block data class (in the example: \texttt{McCormickEnvelopeData}) and the name of the block type (in the example: \texttt{McCormickEnvelope}) do not need to match, but it is a good practice to name the block data class with the same name as the block type and suffix \texttt{Data}.

The methods that we define in the block data class are avaliable in the block type, e.g. Listing~\ref{lst:pyomo.example.2} defines a \texttt{build(x, y)} method that creates the inequalities and adds them to a \texttt{ConstraintList}.

Users can use our custom Pyomo block by first creating a new instance of our block type and then calling any initialization methods we defined, in our example when they call the \texttt{build(x, y)} method the block automatically creates the four inequalities of the McCormick envelope.

\begin{lstlisting}[%
caption={This example shows how to define and then use a custom Pyomo block. The \texttt{McCormickEnvelope} block automatically generates the four inequalities of the McCormick envelope of $w = x y$.},%
captionpos=b,%
label=lst:pyomo.example.2,%
language=Python,%
frame=tb,%
style=ExStyle]
import pyomo.environ as pe
from pyomo.core.base.block import declare_custom_block, _BlockData


@declare_custom_block('McCormickEnvelope')
class McCormickEnvelopeData(_BlockData):
    def __init__(self, component):
        super().__init__(component)
        self.w = None
        self.inequalities = None

    def build(self, x, y):
        """Build envelope of w = x*y. Assume x, y are bounded."""
        del self.inequalities
        del self.w

        self.inequalities = pe.ConstraintList()
        w = pe.Var()
        self.w = w

        xl, xu = x.bounds
        yl, yu = y.bounds

        self.inequalities.add(w >= xl*y + x*yl - xl*yl)
        self.inequalities.add(w >= xu*y + x*yu - xu*yu)
        self.inequalities.add(w <= xu*y + x*yl - xu*yl)
        self.inequalities.add(w <= xl*y + x*yu - xu*yl)


model = pe.ConcreteModel()

model.v1 = pe.Var(bounds=(-2, 2))
model.v2 = pe.Var(bounds=(-3, 1))

# Build envelope for v1*v2
model.envelope = McCormickEnvelope()
model.envelope.build(x=model.v1, y=model.v2)

# Use auxiliary variable w
model.objective = pe.Objective(expr=model.envelope.w)
\end{lstlisting}

\section{Network structure}\label{app:network_structure}

This appendix shows how the pooling library builds and stores a problem. Internally, the network library uses the Python \texttt{networkx} library to store the network. The network is composed of nodes and edges. Node properties:

\begin{enumerate}
\item \texttt{name}: an unique name for the node, used to lookup nodes in the network.
\item \texttt{layer}: nodes are grouped by layer.
\item \texttt{capacity}: a tuple containing the lower and upper capacity of the node (we assume that the bounds on the capacity are non-negative real numbers). If the lower capacity is \texttt{None}, it's assumed to be $0$, if the upper capacity is \texttt{None}, it's assumed to be $\infty$.
\item \texttt{cost}: a cost associated with the node, e.g. the cost of raw materials or the profit for products.
\item \texttt{attr}: a dictionary of additional node attributes. Can be used to store problem-specific properties, e.g. the chemical qualities of a material feed.
\end{enumerate}

Edge properties:

\begin{enumerate}
\item \texttt{source}: the name of the source node.
\item \texttt{destination}: the name of the destination node.
\item \texttt{cost}: the per-unit cost of flow through this edge.
\item \texttt{fixed\_cost}: fixed cost of the edge.
\item \texttt{capacity}: a tuple containing the edge lower and upper capacity. \texttt{None} behaves the same way as node capacity does.
\item \texttt{attr}: a dictionary of additional attributes of the edge.
\end{enumerate}

The network class provides the following methods:

\begin{enumerate}
\item \texttt{add\_node(layer, name, capacity\_lower, capacity\_upper, cost, attr)}: builds a new node and adds it to the network. The node is not connected to any other node yet.
\item \texttt{add\_edge(src, dest, capacity\_lower, capacity\_upper, cost, fixed\_cost, attr)}: adds an edge to the network connecting the \texttt{src} and \texttt{dest} nodes, the nodes must exist in the network or the library will raise an exception.
\item \texttt{successors(src)}: returns an iterator over the nodes that are the successors of \texttt{src}, that is the nodes where exists an edge from \texttt{src} to the node. The method accepts an optional \texttt{layer} argument that filters the nodes that belong to the specified layer.
\item \texttt{predecessors(src)}: returns an iterator over the nodes that are the predecessors of \texttt{src}, that is the nodes where exists an edge from the node to \texttt{src}. The method accepts an optional \texttt{layer} argument that filters the nodes that belong to the specified layer.
\item \texttt{nodes}: returns a dictionary-like object used to access nodes by name.
\item \texttt{edges}: returns a dictionary-like object used to access the edges of the network by source and destination nodes.
\end{enumerate}

As an example, Figure~\ref{fig:pooling.problem.network.library} shows how the library represents the \texttt{adhya4} test problem \citep{adhya1999lagrangian}. Listing~\ref{lst:pooling.problem.network.library} shows how to use the library to build this same problem.

\begin{lstlisting}[%
caption={This code listing shows how to use the pooling network library to create a pooling problem network. Users can add nodes grouped in different ``layers'', for the pooling problem layer 0 contains inputs, layer 1 pools, and layer 2 outputs. Users can add edges between nodes with an optional maximum flow capacity, cost per unit of flow, and fixed cost.},%
captionpos=b,%
label=lst:pooling.problem.network.library,%
language=Python,%
frame=tb,%
style=ExStyle]
import pyomo.environ as pe
from pooling_network.network import Network

# Create network object
network = Network('adhya4')

# Add product inputs.
network.add_node(
    layer=0,
    name='c1',
    capacity_lower=0.0,
    capacity_upper=85.0,
    cost=15.0,
    attr={
        # Quality is a special attribute used by the PoolingPQFormulation block.
        'quality': {
            'q1': 0.5,
            'q2': 1.9,
            'q3': 1.3,
            'q4': 1.0,
        }
    }
)
# Add remaining 7 product inputs c2,...,c8.

# Add pools.
network.add_node(
    layer=1,
    name='o1',
    capacity_lower=0.0,
    capacity_upper=85.0,
    cost=0.0,
    attr=dict(),
)
# Add other pool o2.

# Add outputs.
network.add_node(
    layer=2,
    name='p1',
    capacity_lower=0.0,
    capacity_upper=15.0,
    cost=10.0,
    attr={
        # If quality_upper is present PoolingPQFormulation block will add
        # constraints for the products maximum qualities.
        # If quality_lower is also present, the block will add extra constraints
        # with minimum product qualities.
        'quality_upper': {
            'q1': 1.2,
            'q2': 1.7,
            'q3': 1.4,
            'q4': 1.7,
        }
    }
)
# Add other outputs p2,...,p5.

# Add edges
network.add_edge('c1', 'o1', capacity_upper=85.0, cost=0.0)
# Add all other edges.
\end{lstlisting}

\section{Code}\label{app:code}

\begin{lstlisting}[%
caption={This code listing shows how to use the pooling network library to 1) load the network structure for a pooling problem and then build the \texttt{PoolingPQFormulation} block, 2) add additional constraints that use the flow variables presented in Table~\ref{table:block.formulation}.},%
captionpos=b,%
label=lst:network.add.cuts,%
language=Python,%
frame=tb,%
style=ExStyle]
import pyomo.environ as pe
from pooling_network.instances.literature import literature_problem_data
from pooling_network.instances.data import pooling_problem_from_data
from pooling_network.block import PoolingPQFormulation

# 1) Load adhya4 network structure and build the problem PQ-formulation
problem = pooling_problem_from_data(literature_problem_data('adhya4'))
model = pe.ConcreteModel()
model.pooling = PoolingPQFormulation()
model.pooling.set_pooling_problem(problem)
model.pooling.rebuild()
model.pooling.add_objective(use_flow_cost=False)

# 2) Add additional constraints that use flow variables from the PQ-formulation
b = pe.Block()
##  If there is flow from 'o1' to 'p1', then deactivate 'o1' to 'p2'.
b.has_flow = pe.Var(domain=pe.Binary)
b.set_has_flow = pe.Constraint(expr=(
    model.pooling.y['o1', 'p1'] <= b.has_flow * model.pooling.y['o1', 'p1'].ub
))
b.disable_o1_p2 = pe.Constraint(expr=(
    model.pooling.y['o1', 'p2'] <= (1 - b.has_flow) * model.pooling.y['o1', 'p2'].ub
))
model.additional_cons = b
\end{lstlisting}

\begin{lstlisting}[%
caption={This code listing shows how to use the convenience methods provided by the \texttt{PoolingPQFormulation} block to add valid inequalities and cuts to the pooling problem. We build the linear relaxation of the pooling problem PQ-formulation using \galini{}, but users can use any method to relax the model. The listing shows how, after solving the linear relaxation, the \texttt{PoolingPQFormulation} block can add valid cuts at the solution point.},%
captionpos=b,%
label=lst:pooling.network.library.add.cuts,%
language=Python,%
frame=tb,%
style=ExStyle]
import pyomo.environ as pe
from pooling_network.instances.literature import literature_problem_data
from pooling_network.instances.data import pooling_problem_from_data
from pooling_network.block import PoolingPQFormulation
from galini.relaxations.relax import RelaxationData, relax

# Build Adhya4 network
problem = pooling_problem_from_data(literature_problem_data('adhya4'))

model = pe.ConcreteModel()

# Add a block with the PQ-formulation
model.pooling = PoolingPQFormulation()
model.pooling.set_pooling_problem(problem)
model.pooling.rebuild()

# Add objective function
model.pooling.add_objective(use_flow_cost=False)

lp_solver = pe.SolverFactory('cplex')

# Create a linear relaxation of model using GALINI.
relaxation_data = RelaxationData(model)
linear_model = relax(model, relaxation_data)

# Add variables and inequalities for all
# (pool, output, quality) triplets.
linear_model.pooling.add_inequalities()

for iter in range(10):
    # Solve model
    lp_solver.solve(linear_model)
    print('Iter {}: {}'.format(iter, pe.value(linear_model.pooling.cost)))
    # Automatically add valid cuts for all
    # (pool, output, quality) triplets.
    # This method returns the list of cuts that were added to the problem.
    new_cuts = linear_model.pooling.add_cuts()
    print('  Adding {} cuts'.format(len(new_cuts)))
    if not new_cuts:
        break
\end{lstlisting}

\begin{lstlisting}[%
caption={This code listing shows how to use the Python context manager provided by the network library to add (and automatically remove) the MIP restriction described in Section~\ref{sect:network.primal.heuristic}. Users can specify the number of auxiliary pools each pool is split, and the weight $\gammalt$ of the total flow into pool $l$ each auxiliary pool $t$ receives.},%
captionpos=b,%
label=lst:network.add.heuristic,%
language=Python,%
frame=tb,%
style=ExStyle]
import pyomo.environ as pe
from pooling_network.instances.literature import literature_problem_data
from pooling_network.instances.data import pooling_problem_from_data
from pooling_network.block import PoolingPQFormulation
from pooling_network.heuristic import mip_heuristic

# Build Adhya4 network
problem = pooling_problem_from_data(literature_problem_data('adhya4'))
model = pe.ConcreteModel()
# Add a block with the PQ-formulation
model.pooling = PoolingPQFormulation()
model.pooling.set_pooling_problem(problem)
model.pooling.rebuild()
# Add objective function
model.pooling.add_objective(use_flow_cost=False)
# Add a sub-block that contains the MIP restriction
with mip_heuristic(model.pooling, model.pooling.pooling_problem, tau=2):
    # Print model, it contains a sub-block with the MIP restriction variables
    # and constraints
    model.pprint()
# Automatically removes the sub-block
\end{lstlisting}

\end{document}